\documentclass[11pt]{article}

\usepackage[utf8]{inputenc}
\usepackage[none]{hyphenat} 
\usepackage{amsmath}
\usepackage{amssymb}
\usepackage{amsfonts}
\usepackage{multirow}
\usepackage{enumerate}
\usepackage{hyperref}
\usepackage{graphicx}
\hypersetup{colorlinks=true, linkcolor=blue, urlcolor=blue, citecolor=red}
\usepackage{listings}
\usepackage{color}
\usepackage{rotating}
\usepackage[top=2cm, bottom=2cm, left=2cm, right=2cm]{geometry}
\usepackage{algorithm}
\usepackage{algorithmic}

\newtheorem{example}{Example}
\newtheorem{remark}{Remark}

\title{Estimation of ill-conditioned models using penalized sums of squares of the residuals}
\author{
    \bf{Román Salmerón Gómez}\thanks{Professor, Department of Quantitative methods for economics and business, University of Granada, Spain (e-mail: romansg@ugr.es).}
    \and
    \bf{Catalina B. García García}\thanks{Professor, Department of Quantitative methods for economics and business, University of Granada, Spain (e-mail: cbgarcia@ugr.es).}
}
\date{\today}

\begin{document}

\sloppy

\maketitle

\begin{abstract}
This paper analyzes the estimation of econometric models by penalizing the sum of squares of the residuals with a factor that makes the model estimates approximate those that would be obtained when considering the possible simple regressions between the dependent variable of the econometric model and each of its independent variables. It is shown that the ridge estimator is a particular case of the penalized estimator obtained, which, upon analysis of its main characteristics, presents better properties than the ridge especially in reference to the individual boostrap inference of the coefficients of the model and the numerical stability of the estimates obtained. This improvement is due to the fact that instead of shrinking the estimator towards zero, the estimator shrinks towards the estimates of the coefficients of the simple regressions discussed above.
\end{abstract}


\section{Introduction}

Given the following multiple linear regression model for $n$ observations and $p$ independent variables (constant term included):
\begin{equation}
    \label{model0}
    \mathbf{y}_{n \times 1} = \mathbf{X}_{n \times p} \cdot \boldsymbol{\beta}_{p \times 1} + \mathbf{u}_{n \times 1},
\end{equation}
where $\mathbf{u}$ is the random disturbance (which is assumed to be spherical), the most common estimation method is Ordinary Least Squares (OLS), which consists of minimizing with respect to $\boldsymbol{\beta}$ the sum of squares of the residuals given by the expression:
$$SCR(\boldsymbol{\beta}) = (\mathbf{y} - \mathbf{X} \cdot \boldsymbol{\beta})^{t} \cdot (\mathbf{y} - \mathbf{X} \cdot \boldsymbol{\beta}).$$

However, there are several situations where the objective function to be minimized has been modified, as in the case of the method of Restricted Least Squares (RSM). In this case, the modification seeks to incorporate certain information into the model in order to obtain a more efficient estimation.

It is also common to modify it in order to achieve greater stability in the  estimates of model (\ref{model0}). In this sense, according to O'Driscoll and Ramirez \cite{ODriscollRamirez}, it was Tikhonov \cite{Tikhonov1943} who first proposed the idea of regularization to estimate ill-conditioned problems by introducing additional information.

Other examples with alternative objective functions are the ridge estimator (Hoerl and Kennard \cite{HKa, HKb}) or LASSO regression (Tibshirani \cite{Tibshirani:1996,Tibshirani:2011}), widely used to mitigate the approximate multicollinearity existing in the linear model.

In the former,  $SCR(\boldsymbol{\beta})$ is minimized subject to the constraint $|| \boldsymbol{\beta} ||_{2}^{2} = \sum \limits_{i=1}^{p} \beta_{i}^{2} = \boldsymbol{\beta}^{t} \boldsymbol{\beta} \leq r$ where $r$ is fixed, and implies minimizing with respect to  $\boldsymbol{\beta}$ the Lagrangian:
$$L = SCR(\boldsymbol{\beta}) + k \cdot \boldsymbol{\beta}^{t} \boldsymbol{\beta}, \quad k \geq 0.$$
The solution obtained corresponds to the following expression:
\begin{equation}
    \label{est.cresta}
    \widehat{\boldsymbol{\beta}}(k) = \left( \mathbf{X}^{t} \mathbf{X} + k \cdot \mathbf{I} \right)^{-1} \cdot \mathbf{X}^{t} \mathbf{y},
\end{equation}
where $\mathbf{I}$ denoted the identity matrix of appropriate dimensions. In addition, its variance-covariance matrix is:
\begin{equation}
    var \left( \widehat{\boldsymbol{\beta}}(k) \right) = \sigma^{2} \cdot \left( \mathbf{X}^{t} \mathbf{X} + k \cdot \mathbf{I} \right)^{-1} \cdot \mathbf{X}^{t} \mathbf{X} \cdot \left( \mathbf{X}^{t} \mathbf{X} + k \cdot \mathbf{I} \right)^{-1}.
    \label{est.cresta.var-cor}
\end{equation}

In the second, $SCR(\boldsymbol{\beta})$ is minimized subject to the constraint $|| \boldsymbol{\beta} ||_{1}^{2} = \sum \limits_{i=1}^{p} | \beta_{i} | \leq s$ where $s$ is fix which implies minimizing  with respect to $\boldsymbol{\beta}$ the Lagrangian:
$$L = SCR(\boldsymbol{\beta}) + k \cdot \sum \limits_{i=1}^{p} | \beta_{i} |, \quad k \geq 0.$$

The difference between the two lies in the fact that due to the behavior of the  $l_{1}$ norm,  LASSO contracts and selects variables, while ridge only contracts. It can also be seen that both methods can be viewed as particular situations of the general case in which the penalty added in the objective function is considered to be given by the $l_{q}$ norm, i.e., it is of the form  $|| \boldsymbol{\beta} ||_{q}^{2} = \sum \limits_{i=1}^{p} | \beta_{i} |^{q}$.
Other shrinkage estimators can be found in Liu \cite{Liu1993, Liu2003}, Gao and Liu \cite{GaoLiu2011} or Liu and Gao \cite{LiuGao2011}.

More recently, the methodology proposed by Zou and Hastie \cite{ZouHastie}, known as elastic-net, has emerged as a particularly noteworthy approach. This technique allows the contraction of estimators and the automatic selection of variables by combining the penalty factors (in the sum of squares of the residuals) of the ridge and LASSO regressions. It is also particularly useful when the number of independent variables is much higher than the number of observations and to deal with the problem of multicollinearity (see, for example, Bingzhen et al. \cite{BingzhenWenjuanLingchen2022} or Zou and Zhang \cite{ZouZhang2009}).

Since its publication, this work has been received\footnote{
    Figures as of May 8, 2024 in \href{https://www.scopus.com/record/display.uri?eid=2-s2.0-16244401458&origin=resultslist}{Scopus} and \href{https://scholar.google.es/citations?view_op=view_citation&hl=es&user=tQVe-fAAAAAJ&citation_for_view=tQVe-fAAAAAJ:Tyk-4Ss8FVUC}{Google Scholar}.
} 11820 citations according to SCOPUS and 20777 according to Google Scholar (2260 in 2021, 2345 in 2022 and 2100 in 2023). These figures show that the development of techniques to deal with multicollinearity and variable selection is currently an interesting econometric research topic, probably due to the fact that  \textit{are essential tools in high-dimensional data analysis} (see Ding et al. \cite{DingPengSongChen2023}).

Within this line of research, Wang et al. \cite{Wangetal2021} proposes a novel method called average least squares-centered penalized regression (ALPR). This penalized regression method shrinks the parameters toward the weighted average OLS estimator. This paper provides a theoretical development to obtain the mean square error of the proposed estimator, as well as a simulation study comparing this estimator with the estimators provided by OLS and ridge regression. Some issues to be considered with respect to this paper are the following:
\begin{itemize}
    \item If the OLS estimates have coefficients with different signs, it is possible that the weighted average is close to zero, so there is not much difference between the penalized estimator they propose and the ridge regression.
    \item Shrinking the penalized estimators around the weighted average OLS estimator is just as arbitrary as shrinking them around zero (ridge regression).
    \item ALPR shrinks all the parameters to the same single value. Again, this is a completely arbitrary question without any theoretical justification.
    \item The intercept is not penalized by ALPR, then the approximate non-essential multicollinearity (see footnote in Example \ref{ejemplo.SCRP0} about essential and non-essential multicollinearity) is ignored.
\end{itemize}
Recently, Lukman et al. \cite{Lukman2024} have integrated ALPR with the principal component dimension reduction method to improve its performance in terms of root mean square error.

Alternatively, the present work proposes to minimize the Lagrangian with respect to  $\boldsymbol{\beta}$ the Langragian:
\begin{equation}
    \label{lagran.1}
    L = SCR(\boldsymbol{\beta}) + k \cdot (\boldsymbol{\beta} - \boldsymbol{\alpha})^{t} (\boldsymbol{\beta} - \boldsymbol{\alpha}), \quad k \geq 0,
\end{equation}
where $\boldsymbol{\alpha}$ is a vector that collects the estimates of the slope in all simple linear regressions of the form $\mathbf{y} = \alpha_{1} \cdot \mathbf{1} + \alpha_{i} \cdot \mathbf{X}_{i} + \mathbf{v}$ with $i=2,\dots,p$, where $\mathbf{1}$ denotes the constant term.

    The choice of the parameter around which to shrink the penalized estimators, $\boldsymbol{\alpha}$, is intended to overcome the drawbacks of ALPR highlighted above (see section \ref{rest.SCRP}). In addition, the obtained estimator is analyzed in depth (sections \ref{est.SCRP} to \ref{stability}), providing a consolidated methodology.

Specifically, the work is structured as follows: Section \ref{rest.SCRP} justifies the usefulness of the proposal made to mitigate the approximate multicollinearity existing in the model (\ref{model0}) and achieve stability in the estimation of the coefficients of the independent variables, while in section \ref{est.SCRP} we obtain the expression of the estimator provided by minimizing the lagrangian given in (\ref{lagran.1}) and its main characteristics: trace, norm, variance-covariance matrix, goodness-of-fit and mean square error. To highlight in this section the Monte Carlo simulation performed to analyze the behavior of the mean square error of the proposed penalized estimator.
The implications of this methodology in multicollinearity mitigation are discussed in section \ref{model.aumentado.SCRP}.
The above results are used in section \ref{elec.SCRP} to propose procedures for determining the ideal value of the penalty parameter, $k$. An interpretation of the same is also provided.
From these values, the section \ref{inf.bootstrap} shows how to perform inference from bootstrap methods.
Section \ref{stability} analyzes the stability of the estimates obtained from the proposed penalized estimator and section \ref{ejemplos.SCRP} illustrates the results obtained by means of an example with real data.
The paper ends with the section \ref{conc.SCRP}, where the main results obtained are highlighted. Finally, the code utilized in R \cite{RCoreTeam} to generate the results presented in this study is accessible on GitHub, specifically at \url{https://github.com/rnoremlas/Penalized-estimator}.

\section{Justification of the proposed restriction} \label{rest.SCRP}

Taking into account the following notions of partial and full effect given in Novales \cite{NovalesWeb}:
\begin{itemize}
    \item \textit{What is the impact on $\mathbf{y}$ of a unit variation in $\mathbf{X}_{i}$ if the other explanatory variables in the model did not vary? Answer: partial effect (multiple linear regression).}
    \item \textit{What is the total impact on $\mathbf{y}$ of a unit variation in $\mathbf{X}_{i}$ if the rest of the explanatory variables vary as would be expected given the correlations observed between them throughout the sample? Answer: total effect (simple linear regression).}
\end{itemize}

Both effects should coincide (see example \ref{ejemplo.SCRP0}) if the \textit{ceteris paribus}, i.e., when one variable varies the rest remain constant, were really verified.

\begin{example}
    \label{ejemplo.SCRP0}

    Given the model $\mathbf{y} = \beta_{1} \cdot \mathbf{1} + \beta_{2} \cdot \mathbf{X}_{2} + \beta_{3} \cdot \mathbf{X}_{3} + \mathbf{u}$, 100 observations are simulated assuming that $\mathbf{X}_{2}$ and $\mathbf{X}_{3}$  are distributed according to a normal distribution of mean 1 and variance 100, $\mathbf{u}$ according to a normal distribution of mean 0 and variance 2 and the dependent variable $\mathbf{y}$ is generated as $\mathbf{y} = 5 \cdot \mathbf{1} + 2 \cdot \mathbf{X}_{2} - 4 \cdot \mathbf{X}_{3} + \mathbf{u}$.

    It is verified that the coefficient of simple correlation between $\mathbf{X}_{2}$ and $\mathbf{X}_{3}$ is equal to -0.01684487. This leads to a variance inflation factor equal to 1.000284 which is very close to its minimum value of 1. This indicates that the degree of approximate multicollinearity of the essential type\footnote{
        Marquardt \cite{Marquardt1980}, Marquardt and Snee \cite{MarquardtSnee} and Snee and Marquardt \cite{SneeMarquardt} distinguish between approximate multicollinearity of the essential type (relationship between the independent variables of the model excluding the constant term) and non-essential (relationship between the constant term and at least one of the remaining independent variables of the model). Salmerón et al. \cite{Salmeron2020} shows that the variance inflation factor only detects variance of the essential type and that the coefficient of variation must be used to detect the non-essential multicollinearity.
         According to these authors, a coefficient of variation  below 0.1002506 indicates that  the existing approximate multicollinearity of nonessential type is of concern.} is not troubling. The high values for the coefficients of variation,  5.898071 for $\mathbf{X}_{2}$ and 3.022512 for $\mathbf{X}_{3}$, indicate that the nonessential type is not troubling either.

    OLS estimation provides the estimates $\widehat{\beta}_{1} = 4.658938$, $\widehat{\beta}_{2} = 1.998514$ and $\widehat{\beta}_{3} = -3.99332$. Note that the estimates practically coincide with the real coefficients.

    On the other hand, given the simple linear regressions $\mathbf{y} = \alpha_{1} \cdot \mathbf{1} + \alpha_{i} \mathbf{X}_{i} + \mathbf{v}$, with $i=2,3$, the following OLS estimators are obtained: $\widehat{\alpha}_{2} = 2.061951$ and $\widehat{\alpha}_{3} = -4.029017$. Once again, very similar to the true values of $\beta_{2}$ and $\beta_{3}$.
    \hfill $\Box$
\end{example}

Developing the previous example from a theoretical point of view:
\begin{itemize}
    \item  Estimates of simple linear regressions $\mathbf{y} = \alpha_{1} \cdot \mathbf{1} + \alpha_{i} \cdot \mathbf{X}_{i} + \mathbf{v}$ with $i=2,\dots,p$ are given by the expression:
        $$\widehat{\alpha}_{i} = \frac{cov(\mathbf{y}, \mathbf{X}_{i})}{var (\mathbf{X}_{i})}, \quad i=2,\dots,p.$$
    \item Assuming orthogonality in the model  (\ref{model0}),  the OLS estimator responds to the expression:
        \begin{eqnarray*}
            \widehat{\boldsymbol{\beta}} &=& \left( \mathbf{X}^{t} \mathbf{X} \right)^{-1} \mathbf{X} \mathbf{y} \\
                &=& \left( \begin{array}{cccc}
                    n & 0 & \cdots & 0 \\
                    0 & \sum \limits_{j=1}^{n} X_{2j}^{2} & \cdots & 0 \\
                    \vdots & \vdots & \ddots & \vdots \\
                    0 & 0 & \cdots & \sum \limits_{j=1}^{n} X_{kj}^{2} \\
                \end{array} \right)^{-1} \left(
                \begin{array}{c}
                    \sum \limits_{j=1}^{n} y_{j} \\
                    \sum \limits_{j=1}^{n} X_{2j} y_{j} \\
                    \vdots \\
                    \sum \limits_{j=1}^{n} X_{kj} y_{j} \\
                \end{array} \right)^{-1} = \left(
                \begin{array}{c}
                    \overline{\mathbf{y}} \\
                    \frac{cov(\mathbf{y}, \mathbf{X}_{2})}{var (\mathbf{X}_{2})} \\
                    \vdots \\
                    \frac{cov(\mathbf{y}, \mathbf{X}_{k})}{var (\mathbf{X}_{k})} \\
                \end{array} \right) \\
        \end{eqnarray*}
        since, for $i=2,\dots,p$, it is verified that $cov(\mathbf{y}, \mathbf{X}_{i}) = \frac{1}{n} \sum \limits_{j=1}^{n} X_{ij} y_{j}$ and $var (\mathbf{X}_{i}) = \frac{1}{n} \sum \limits_{j=1}^{n} X_{ij}^{2}$  finding that $\sum \limits_{j=1}^{n} X_{2j} = \dots = \sum \limits_{j=1}^{n} X_{pj} = 0$ by the orthogonality condition between  $\mathbf{X}_{i}$ and the intercept.
\end{itemize}

For this reason, in the expression (\ref{lagran.1}) it will be considered that:
\begin{equation}
    \boldsymbol{\alpha}_{p \times 1} = \left( \overline{\mathbf{y}}, \ \frac{cov(\mathbf{y}, \mathbf{X}_{2})}{var (\mathbf{X}_{2})}, \dots, \ \frac{cov(\mathbf{y}, \mathbf{X}_{p})}{var (\mathbf{X}_{p})} \right)^{t},
    \label{alfa}
\end{equation}
that is, $\boldsymbol{\alpha}$ is a known fixed vector.

In short, in the restriction imposed on the sum of squares of the penalized residuals represented in the lagrangian of the expression (\ref{lagran.1}) it will be considered that $\boldsymbol{\alpha}$ will take the value given in (\ref{alfa}). In this way, the purpose is to make the estimates of the coefficients of the multiple linear regression model as close as possible to the estimates obtained in the corresponding simple linear regressions.

\begin{example}[continuation of Example \ref{ejemplo.SCRP0}]
    \label{ejemplo.SCRP}

    If, on the other hand, it is generated $\mathbf{X}_{3}$ como $\mathbf{X}_{3} = \mathbf{1} + 5 \cdot \mathbf{X}_{2} + \widetilde{\mathbf{1}}$, where $\widetilde{\mathbf{1}}_{i} = (-1)^{i},$ it is verified that the coefficient of simple correlation between $\mathbf{X}_{2}$ and $\mathbf{X}_{3}$ is equal to 0.9997975, which leads to a variance inflation factor equal to 2468.787. In this case, the coefficient of variation of $\mathbf{X}_{3}$ is equal to 5.279348. Thus, the approximate multicollinearity of the essential type is troubling while the non-essential type is not.

    The estimates obtained in this case are as follows:
    $$\widehat{\beta}_{1} = 4.3740148, \quad \widehat{\beta}_{2} = 0.4443309, \quad \widehat{\beta}_{3} = -3.6897307, \quad \widehat{\alpha}_{2} = -18.03679, \quad \widehat{\alpha}_{3} = -3.601057,$$
    being the estimates of $\beta_{1}$ and $\beta_{3}$ significantly different from zero (with 5\% of significance level), contrary to that of $\beta_{2}$.

    It is found that the estimation of $\widehat{\beta}_{2}$ differs from the true value of the coefficient of $\mathbf{X}_{2}$, while the value of $\widehat{\alpha}_{2}$ captures the true value of the coefficient obtained by substituting $\mathbf{X}_{3} = \mathbf{1} + 5 \cdot \mathbf{X}_{2} + \widetilde{\mathbf{1}}$ en $\mathbf{y} = 5 \cdot \mathbf{1} + 2 \cdot \mathbf{X}_{2} - 4 \cdot \mathbf{X}_{3} + \mathbf{u}$. It is to say, $\mathbf{y} = \left( \mathbf{1} - 4 \cdot \widetilde{\mathbf{1}} \right) - 18 \cdot \mathbf{X}_{2} + \mathbf{u}$: as captured by simple regression rather than multiple regression, the real relationship between $\mathbf{X}_{2}$ and  $\mathbf{y}$ is actually -18 rather than 2.

    To summarize, it can be observed that in the first case, the degree of essential approximate multicollinearity is not troubling while in the second case it is. Note also that the consequence of high linear relationships between the independent variables is to obtain different values for the partial and total effect.
    \hfill $\Box$
\end{example}

\section{Penalized estimator} \label{est.SCRP}

In this section we obtain the estimator that verifies the above conditions, as well as its norm, variance-covariance matrix, goodness-of-fit and mean square error. We also consider the augmented model whose OLS estimation leads to the penalized estimator obtained. In the latter case, the effect of this augmented model on the detection of the degree of multicollinearity existing once the proposed estimator has been applied is analyzed.

With this purpose, instead of working with the expression (\ref{lagran.1}), expression (\ref{lagran.2}) will be minimized where previously a parameter $h$ has been introduced as follows:
\begin{equation}
    \label{lagran.2}
    L = SCR(\boldsymbol{\beta}) + k \cdot (\boldsymbol{\beta} - h \cdot \boldsymbol{\alpha})^{t} (\boldsymbol{\beta} - h \cdot \boldsymbol{\alpha}), \quad k \geq 0,
\end{equation}
note that for $h=0$ the Hoerl and Kennard estimator would be obtained. \cite{HKa, HKb}. That is, the ridge estimator is a particular case of the one presented here.

While it makes sense that $h \in [0, \ 1]$, it will be considered only that $h \in \{0, 1\}$ (for further details see subsection \ref{h_0_1}).

\subsection{Estimation}

Minimizing with respect to $\beta$ the expression (\ref{lagran.2}) the following estimator is obtained:
\begin{equation}
    \widehat{\boldsymbol{\beta}}(k, h) = \left( \mathbf{X}^{t} \mathbf{X} + k \cdot \mathbf{I} \right)^{-1} \cdot \left( \mathbf{X}^{t} \mathbf{y} + k \cdot h \cdot \boldsymbol{\alpha} \right), \label{est.SCRP.1}
\end{equation}
where it is clear that if $h=0$, the penalized estimator coincides with the ridge (the expressions (\ref{est.cresta}) and (\ref{est.SCRP.1}) are the same) and if $k=0$ the OLS estimation is obtained.

\begin{remark}
    Note that if it is verified that $\boldsymbol{\alpha} = \mathbf{0}$, where $\mathbf{0}$ is a vector of zeros of appropriate dimensions, the penalized estimator also coincides with the ridge estimator. That is, the ridge estimator could be considered to be the particular case of the penalized estimator in which all the coefficients of the possible simple regressions are zero and the dependent variable is centered:
    $$\overline{\mathbf{y}} = 0, \  cov(\mathbf{y}, \mathbf{X}_{i}) = 0 \mbox{ for } i=2,\dots,p.$$

  This implies that the total impact of the variations of each independent variable (excluding the constant term) on the dependent variable is null, assuming that the rest of the independent variables vary as would be expected given the correlations observed between them in the sample.
\end{remark}

\begin{remark}

        Note that in the case where $h=1$ and $\boldsymbol{\alpha} = \mathbf{d} \left( \mathbf{d}^{t} \mathbf{X}^{t} \mathbf{X} \mathbf{d} \right)^{-1} \mathbf{d}^{t} \mathbf{X}^{t} \mathbf{X} \widehat{\boldsymbol{\beta}}$ where $\mathbf{d}$ is a p-dimensional vector in which all the elements are one and $\widehat{\boldsymbol{\beta}}$ is the OLS estimator, the ALPR estimator is obtained.

        Similarly, for $h=1$, $d=-k$ and $\boldsymbol{\alpha} = \widehat{\boldsymbol{\beta}}$, the expression (\ref{est.SCRP.1}) coincides with the Liu \cite{Liu2003} estimator.

\end{remark}

On the other hand, denoting $\mathbf{Z} (k) = \left( \mathbf{X}^{t} \mathbf{X} + k \cdot \mathbf{I} \right)^{-1}$ and taking into account that $\mathbf{X}^{t} \mathbf{X} \cdot \widehat{\boldsymbol{\beta}} = \mathbf{X}^{t} \mathbf{y}$, it is obtained that expression (\ref{est.SCRP.1}) can be rewritten as:
\begin{equation}
    \widehat{\boldsymbol{\beta}}(k, h) = \mathbf{Z} (k) \cdot \mathbf{X}^{t} \mathbf{X} \cdot \widehat{\boldsymbol{\beta}} + k \cdot h \cdot \mathbf{Z} (k) \cdot \boldsymbol{\alpha}, \label{est.SCRP.2}
\end{equation}
and in this case:
$$E \left[ \widehat{\boldsymbol{\beta}}(k, h) \right] = \mathbf{Z} (k) \cdot \mathbf{X}^{t} \mathbf{X} \cdot \boldsymbol{\beta} + k \cdot h \cdot \mathbf{Z} (k) \cdot \boldsymbol{\alpha},$$
i.e., the estimator given in (\ref{est.SCRP.1}) is a biased estimator of $\boldsymbol{\beta}$ as long as $k \not= 0$.

Finally, note that from the expressions (\ref{est.cresta}) and (\ref{est.SCRP.1}) it is obtained that:
\begin{equation}
    \label{penalizado.cresta}
    \widehat{\boldsymbol{\beta}}(k, h) = \widehat{\boldsymbol{\beta}}(k) + k \cdot h \cdot \mathbf{Z} (k) \cdot \boldsymbol{\alpha}.
\end{equation}

\begin{remark}
    \label{kas}
    Considering that the expression (\ref{lagran.2}) is equivalent to:
    \begin{equation}
        \label{lagran.3}
        L = k_{1} \cdot SCR(\boldsymbol{\beta}) + k_{2} \cdot (\boldsymbol{\beta} - h \cdot \boldsymbol{\alpha})^{t} (\boldsymbol{\beta} - h \cdot \boldsymbol{\alpha}),
    \end{equation}
    just considering that $k = \frac{k_{2}}{k_{1}}$ with $k_{1}>0$.
    Minimize with respect to $beta$ the expression (\ref{lagran.3}) leads to :
    \begin{eqnarray}
        \widehat{\boldsymbol{\beta}}(k_{1}, k_{2},h) &=& \left( k_{1} \cdot \mathbf{X}^{t} \mathbf{X} + k_{2} \cdot \mathbf{I} \right)^{-1} \cdot \left( k_{1} \cdot \mathbf{X}^{t} \mathbf{y} + k_{2} \cdot h \cdot \boldsymbol{\alpha} \right) \label{est.SCRP.1.bis} \\
        & = & \left( k_{1} \cdot \mathbf{X}^{t} \mathbf{X} + k_{2} \cdot \mathbf{I} \right)^{-1} \cdot \left[ \left( k_{1} \cdot \mathbf{I} + \frac{1}{n} \cdot k_{2} \cdot h \cdot \mathbf{V} \right) \cdot \mathbf{X}^{t} \mathbf{y} - k_{2} \cdot h \cdot \overline{\mathbf{y}} \cdot \mathbf{V} \cdot \mathbf{B} \right], \nonumber
    \end{eqnarray}
    donde:
    $$\mathbf{V}_{p \times p} = diag \left( 1, \ \frac{1}{var(\mathbf{X}_{2})}, \dots, \ \frac{1}{var(\mathbf{X}_{p})} \right), \quad
    \mathbf{B}_{k \times 1} = \left( 0, \ \overline{\mathbf{X}}_{2}, \dots, \ \overline{\mathbf{X}}_{p} \right).$$

    Note that if the variables are centered, i.e., $\mathbf{B} = \mathbf{0}_{p \times 1}$ where $\mathbf{0}$ s a vector of zeros of appropriate dimensions, the expression (\ref{est.SCRP.1.bis}) is rewritten as:
    \begin{eqnarray}
        \widehat{\boldsymbol{\beta}}(k_{1}, k_{2},h) &=& \left( k_{1} \cdot \mathbf{X}^{t} \mathbf{X} + k_{2} \cdot \mathbf{I} \right)^{-1} \cdot \left( k_{1} \cdot \mathbf{I} + \frac{1}{n} \cdot k_{2} \cdot h \cdot \mathbf{V} \right) \cdot \mathbf{X}^{t} \mathbf{y}, \nonumber \\
            &=& \left( \mathbf{X}^{t} \mathbf{X} + k \cdot \mathbf{I} \right)^{-1} \cdot \left( \mathbf{I} + \frac{1}{n} \cdot k \cdot h \cdot \mathbf{V} \right) \cdot \mathbf{X}^{t} \mathbf{y}, \label{est.SCRP.2.bis}
    \end{eqnarray}
    where now:
    $$\mathbf{V}_{p \times p} = diag \left( 1, \ \frac{n}{\sum \limits_{i=1}^{n} X_{i2}^{2}}, \dots, \ \frac{n}{\sum \limits_{i=1}^{n} X_{ip}^{2}} \right).$$

    Given a problem of bad conditioning of the matrix $\mathbf{X}$, it will affect the calculation of $\widehat{\boldsymbol{\beta}}$ both in the expression $\left( \mathbf{X}^{t} \mathbf{X} \right)^{-1}$ and in $\mathbf{X}^{t} \mathbf{y}$. In this case, from the expression (\ref{est.SCRP.2.bis}), we would be considering both terms and not only the first one as, for example, the ridge estimator does (see expression (\ref{est.cresta})).
\end{remark}

\subsubsection{Orthogonal case}
    \label{h_0_1}

Section \ref{rest.SCRP} shows that in the orthogonal case it is verified that $\widehat{\boldsymbol{\beta}} = \boldsymbol{\alpha}$, therefore, taking into account that $\mathbf{X}^{t}\mathbf{X} \cdot \widehat{\boldsymbol{\beta}} = \mathbf{X}^{t}\mathbf{y}$ the expression (\ref{est.SCRP.1}) can be rewritten as:
\begin{eqnarray*}
    \widehat{\boldsymbol{\beta}}(k, h) &=& \left( \mathbf{X}^{t}\mathbf{X} + k \cdot \mathbf{I} \right)^{-1} \cdot \left( \mathbf{X}^{t}\mathbf{X} + k \cdot h \cdot \mathbf{I} \right) \cdot \widehat{\boldsymbol{\beta}} \\
    &=& \left( \begin{array}{cccc}
        \frac{x_{11} + k \cdot h}{x_{11} + k} & 0 & \cdots & 0 \\
        0 & \frac{x_{22} + k \cdot h}{x_{22} + k} & \cdots & 0 \\
        \vdots & \vdots & \ddots & \vdots \\
        0 & 0 & \cdots & \frac{x_{pp} + k \cdot h}{x_{pp} + k} \\
    \end{array} \right) \cdot \widehat{\boldsymbol{\beta}},
\end{eqnarray*}
where $x_{ii}$ is denoted as the element (i,i) of $\mathbf{X}^{t}\mathbf{X}$, it is to say, $x_{ii} = \sum \limits_{j=1}^{n} X_{ij}^{2}$ for $i=1,\dots,p$. In this case:
\begin{itemize}
    \item If $h=0$, the penalized estimator coincides with the ridge and it is found that $\widehat{\boldsymbol{\beta}}(k, h) \not= \widehat{\boldsymbol{\beta}}$ for values of $k$ different from zero.
    \item If $h=1$, then $\widehat{\boldsymbol{\beta}}(k, h) = \widehat{\boldsymbol{\beta}}$ for any value of $k$.
    \item If $h \in (0,1)$, then $\widehat{\boldsymbol{\beta}}(k, h) \not= \widehat{\boldsymbol{\beta}}$ for values of $k$ different from zero.
\end{itemize}

Thus, although technically $h$ can be any non-negative value, the only advisable value for $h$ would be 1, since in this case it is assured that the penalized estimator coincides with the OLS estimator when there is orthogonality in the model (\ref{model0}). The case $h=0$ will be also considered to work with the ridge estimator as a particular case of the penalized estimator.

\subsection{Trace and norm of the penalized estimator}

From expression (\ref{est.SCRP.1}) the trace of the estimator can be obtained by simply plotting $\widehat{\boldsymbol{\beta}}(k, h)$ for different values of $k$ (it will be considered that $k \in [0, 1]$). In this way, it is possible to obtain graphically which values of $k$ stabilize the estimates of $\boldsymbol{\beta}$ obtained from $\widehat{\boldsymbol{\beta}}(k, h)$.

Since  section \ref{norma.penalizada} of the appendix \ref{apen.operaciones} shows that:
$$|| \widehat{\boldsymbol{\beta}}(k, h) || \rightarrow h^{2} \cdot || \boldsymbol{\alpha}||,$$
it is guaranteed that there are $k$ values that will stabilize the values provided by the penalized estimator. In addition, as stated in the expression (\ref{lagran.2}), these estimates converge to the values of $\boldsymbol{\alpha}$.

\subsection{Matrix of variances-covariances}

From expression (\ref{est.SCRP.2}) and taking into account that $var \left( \widehat{\boldsymbol{\beta}} \right) = \sigma^{2} \cdot \left( \mathbf{X}^{t} \mathbf{X} \right)^{-1}$ it is immediately clear that:
\begin{equation}
    \label{var.cov.est.SCRP}
    var \left( \widehat{\boldsymbol{\beta}}(k, h) \right) = \sigma^{2} \cdot \mathbf{Z} (k) \cdot \mathbf{X}^{t} \mathbf{X} \cdot \mathbf{Z} (k).
\end{equation}
It is observed that this expression does not depend on the $h$ factor and coincides with that of the Hoerl and Kennard ridge estimator \cite{HKa, HKb} given in expression (\ref{est.cresta.var-cor}).

Section \ref{varianza.optima} of appendix \ref{apen.operaciones} shows that $var \left( \widehat{\boldsymbol{\beta}}(k, h) \right) < var \left( \widehat{\boldsymbol{\beta}} \right)$, i.e., the penalized estimator is optimal in the sense that it has a smaller variance-covariance matrix than the OLS estimator.

\subsection{Goodness-of-fit}

Taking into account that from the expression (\ref{est.SCRP.1}) the residuals of the penalized estimate would be obtained from $\mathbf{e}(k, h) = \mathbf{y} - \mathbf{X} \widehat{\boldsymbol{\beta}}(k, h)$, it is clear that $\mathbf{y} = \mathbf{X} \widehat{\boldsymbol{\beta}}(k, h) + \mathbf{e} (k, h)$. In this case the following equality is verified:
$$\mathbf{y}^{t} \mathbf{y} = \widehat{\boldsymbol{\beta}}(k, h)^{t} \left( \mathbf{X}^{t} \mathbf{X} + 2 \cdot k \cdot \mathbf{I} \right) \widehat{\boldsymbol{\beta}}(k, h) - 2 \cdot k \cdot h \cdot \widehat{\boldsymbol{\beta}}(k, h)^{t} \boldsymbol{\alpha} + \mathbf{e}(k, h)^{t} \mathbf{e}(k, h).$$

Consequently, a measure of goodness-of-fit can be given by:
\begin{eqnarray}
    GoF(k, h) &=& 1 - \frac{\mathbf{e} (k, h)^{t} \cdot \mathbf{e} (k, h)}{\mathbf{y}^{t} \mathbf{y}} \nonumber \\
        &=& \frac{\widehat{\boldsymbol{\beta}}(k, h)^{t} \left( \mathbf{X}^{t} \mathbf{X} + 2 \cdot k \cdot \mathbf{I} \right) \widehat{\boldsymbol{\beta}}(k, h) - 2 \cdot k \cdot h \cdot \widehat{\boldsymbol{\beta}}(k, h)^{t} \boldsymbol{\alpha}}{\mathbf{y}^{t} \mathbf{y}}, \label{ba.ajuste.alterna}
\end{eqnarray}
where $\mathbf{y}^{t} \mathbf{y} = n \cdot var(\mathbf{y})$. Note that
\begin{itemize}
    \item For $k=0$, it is obtained the coefficient of determination obtained, $R^{2}$, by OLS in model (\ref{model0}), as long as the dependent variable has zero mean (see Salmerón et al. \cite{Salmeron2020b} for details).
    \item For $h=0$ there would be a goodness of fit for the ridge estimator.
    \item Since residuals do not necessarily sum to zero (see section \ref{ba.errores} of Appendix \ref{apen.operaciones}),  this measure does not have to vary between zero and one as it does in OLS..
    \item Section \ref{ba.monotonia} of Appendix \ref{apen.operaciones} Analyzes the monotony of $GoF(k, h)$, concluding that:
        \begin{itemize}
            \item if $h \not=0$, its decrease is not assured when $k$ increases, although this is to be expected. In this case, it would be verified that $GoF(k,h) \in (2 \cdot h \cdot \mathbf{y}^{t} \mathbf{X} \boldsymbol{\alpha}, R^{2}]$.
            \item if $h =0$, its decreasing is assured when $k$ increases and it would be verified that $GoF(k,h) \in (0, R^{2}]$.
        \end{itemize}
\end{itemize}

\subsection{Mean Square Error}
    \label{ecm_kh}

Since the penalized estimator is biased, it is interesting to calculate its mean square error (MSE) in order to establish when it is lower than that obtained by OLS. In this case, the MSE corresponds to the following expression:
\begin{equation}
    ECM \left( \widehat{\boldsymbol{\beta}}(k, h) \right) = \sigma^{2} \cdot traza \left( var \left( \widehat{\boldsymbol{\beta}}(k, h) \right) \right) + \left( E \left[ \widehat{\boldsymbol{\beta}}(k, h) \right] - \boldsymbol{\beta} \right)^{t} \left( E \left[ \widehat{\boldsymbol{\beta}}(k, h) \right] - \boldsymbol{\beta} \right), \label{ECM.SCRP.1}
\end{equation}
where $E \left[ \widehat{\boldsymbol{\beta}}(k, h) \right] - \boldsymbol{\beta} = \left( \mathbf{Z}(k) \cdot \mathbf{X}^{t} \mathbf{X} - \mathbf{I} \right) \cdot \boldsymbol{\beta} + k \cdot h \cdot \mathbf{Z} (k) \cdot \boldsymbol{\alpha}$.

Since the variance-covariance matrix of the penalized estimator coincides with that of the ridge estimator, the expression (\ref{ECM.SCRP.1}) can be rewritten as:
\begin{equation}
    ECM \left( \widehat{\boldsymbol{\beta}}(k, h) \right) = ECM \left( \widehat{\boldsymbol{\beta}}(k) \right) + \mathbf{S} (k, h), \label{ECM.SCRP.2}
\end{equation}
where $ECM \left( \widehat{\boldsymbol{\beta}}(k) \right)$ is the MSE of ridge estimator and $\mathbf{S} (k, h) = 2 \cdot k \cdot h \cdot \boldsymbol{\beta}^{t} \cdot \left( \mathbf{Z} (k) \cdot \mathbf{X}^{t} \mathbf{X} - \mathbf{I} \right)^{t} \cdot \mathbf{Z} (k) \cdot \boldsymbol{\alpha} + k^{2} \cdot h^{2} \cdot \boldsymbol{\alpha}^{t} \cdot \mathbf{Z} (k)^{t} \cdot \mathbf{Z} (k) \cdot \boldsymbol{\alpha}$. 

Hoerl and Kennard \cite{HKa} showed that the mean square error of the ridge estimator is equal to the sum of two terms::
$$ECM \left( \widehat{\boldsymbol{\beta}}(k) \right) = \sigma^{2} \sum \limits_{i=1}^{p} \frac{\lambda_{i}}{(\lambda_{i}+k)^{2}} + k^{2} \boldsymbol{\beta}^{t} \left( \mathbf{X}^{t} \mathbf{X} + k \cdot \mathbf{I} \right)^{-2} \boldsymbol{\beta},$$
such that the first is decreasing in $k$ and the second is increasing, therefore, the monotonicity of$ECM \left( \widehat{\boldsymbol{\beta}}(k) \right)$ is not clear. However, due to
$$k^{2} \boldsymbol{\beta}^{t} \left( \mathbf{X}^{t} \mathbf{X} + k \cdot \mathbf{I} \right)^{-2} \boldsymbol{\beta} \longrightarrow \boldsymbol{\beta}^{t}\boldsymbol{\beta},$$
when $k \rightarrow +\infty$, it is verified that $ECM \left( \widehat{\boldsymbol{\beta}}(k) \right)$ has a horizontal asymptote.

Also, section \ref{ba.esgoECM} of Appendix \ref{apen.operaciones} shows that the first term of $\mathbf{S} (k, h)$ has no definite monotonicity while the second is increasing in $k$. Since both terms have a horizontal asymptote, there must be a value of $k$ for which the value of $\mathbf{S}(k, h)$ is stable. This is an important consideration because if it is verified that $\mathbf{S} (k, h) < 0$, then $ECM \left( \widehat{\boldsymbol{\beta}}(k, h) \right) < ECM \left( \widehat{\boldsymbol{\beta}}(k) \right)$.

Finally, it should be noted that $ECM \left( \widehat{\boldsymbol{\beta}}(k, h) \right)$ has a horizontal asymptote because when the $k \rightarrow +\infty$ it is verified\footnote{
    Note that for $h=1$: $ECM \left( \widehat{\boldsymbol{\beta}}(k, h) \right) \longrightarrow \left( \boldsymbol{\beta} - \boldsymbol{\alpha} \right)^{t} \left( \boldsymbol{\beta} - \boldsymbol{\alpha} \right)$ when $k \rightarrow +\infty$.
} that:
$$ECM \left( \widehat{\boldsymbol{\beta}}(k, h) \right) \longrightarrow \boldsymbol{\beta}^{t}\boldsymbol{\beta} - 2 \cdot h \cdot \boldsymbol{\beta}^{t} \boldsymbol{\alpha} + h^{2} \cdot \boldsymbol{\alpha}^{t} \boldsymbol{\alpha}.$$
Therefore, there must be a value of $k$ at which the $ECM \left( \widehat{\boldsymbol{\beta}}(k, h) \right)$ is stabilized.

\subsubsection{Monte Carlo Simulation}
    \label{simulatioMC}

In order to obtain information on the comparison of the ECM of the penalized estimator with that of the OLS and ridge estimator, and due to the impossibility of analytically establishing the behavior of the ECM monotonicity, a Monte Carlo simulation is proposed below where values are simulated from:

$$\mathbf{X}_{i} = \sqrt{1 - \xi^{2}} \cdot \mathbf{W}_{i} + \mathbf{W}_{p}, \quad i=2,\dots,p,$$
where $p=3, 4, 5, 6$, $\xi \in \{ 0.96, 0.97, 0.98, 0.99 \}$, $\mathbf{W}_{i} \sim N(\mu, \sigma^{2})$ con $\mu$ randomly obtained from the set $\{ 0, \pm 2, \pm 4, \pm 6, \pm 8, \pm 10 \}$ and $\sigma \in \{ 0.01, 0.1, 5, 10, 15 \}$ and $n \in \{ 30, 40, 50, \dots, 200 \}$.
This way of simulating data has been previously used, for example, by McDonald and Galarneau \cite{McDonaldGalarneau1975}, Wichern and Churchill \cite{WichernChurchill1978}, Gibbons \cite{Gibbons1981}, Kibria \cite{Kibria2003}, Salmerón et al. \cite{Salmeronetal2016, salmeron2018JSCS} or Rodríguez et al. \cite{Rodriguezetal2021, Rodriguez2022}; and the goal is that each pair of independent variables has a correlation coefficient equal to $\xi^{2}$.

Once the data has been simulated (1440 simulations are performed), the matrix $\mathbf{X} = [ \mathbf{1}, \mathbf{X}_{2}, \dots, \mathbf{X}_{p}]$ is generated and also $\mathbf{y} = \mathbf{X} \boldsymbol{\beta} + \mathbf{u}$ where the elements of $\boldsymbol{\beta}$ are obtained randomly according to $\beta_{i} \in \{ \pm 1, \pm 2, \pm 3, \pm 4, \pm 5 \}$ and $\mathbf{u} \sim N(0, 1)$.

For each case, considering that $k \in \{ 0,0.01, 0.02, \dots, 0.09, 1\}$, it is calculated:
\begin{itemize}
     \item The mean square error in OLS, $ECM_{OLS} = ECM \left( \widehat{\boldsymbol{\beta}}(0, 0) \right)$, the minimum value of the mean square error for the ridge estimator, $ECM_{R} = ECM \left( \widehat{\boldsymbol{\beta}}(k, 0) \right)$ and for the penalized estimator for $h=1$, $ECM_{P} = ECM \left( \widehat{\boldsymbol{\beta}}(k, 1) \right)$. The Algorithm \ref{algoritmo} has been applied in this case. Note that this algorithm also calculates whether the minimum obtained is a unique minimum or not: in all cases where the minimum has been reached, it is unique.
     \item The minimum coefficient of variation (CV), the maximum variance inflation factor (VIF) and the condition number (CN). In this way, the existence of approximate multicollinearity of essential and non-essential type in the simulated data is analyzed. Calculations were performed with the package ``multiColl'' \cite{multiColl} of R \cite{RCoreTeam} (for further details see Salmerón et al. \cite{Salmeron2021,Salmeron2022}).
\end{itemize}

Table \ref{detec_multicol} shows the minimum, mean, median and maximum values for the minimum coefficient of variation, maximum variance inflation factor and condition number. Considering the traditional thresholds for these multicollinearity detection measures, it is observed that situations have been simulated where the existing approximate multicollinearity (both essential and nonessential) is and is not troubling. In addition, the significant difference between the mean and median of the values obtained for the number of condition is notable, indicating the existence of very high values that shift the distribution to the right.

\begin{table}
    \centering
    \begin{tabular}{cccc}
        \hline
        & Minimum CV & Maximum VIF & CN \\
        \hline
        Minimum & 0.00123 & 6.394 & 5.966 \\
        Mean & 2.48557 & 71.689 & 1752.275 \\
        Median & 1.03504 & 52.495 & 33.489 \\
        Maximum & 42.55852 & 367.690 & 25384.87 \\
        \hline
    \end{tabular}
    \caption{Minimum, mean, median and maximum value for the minimum coefficient of variation (CV), the maximum variance inflation factor (VIF) and condition number (CN)} \label{detec_multicol}
\end{table}

\begin{table}
    \centering
    \begin{tabular}{ccccc}
        \hline
        Case & Number (\%) & Minimum CV & Maximum VIF & CN \\
        \hline
        (A) $MSE_{P} < MSE_{R} < MSE_{OLS}$ & 0 (0\%) &  &  & \\
        \cline{3-5}
         &  & 1.85 & 31.199 & 12.726 \\ 
        (B) $MSE_{P} < MSE_{OLS} < MSE_{R}$ & 408 (46.52\%) & 4.6786 & 65.957 & 18.934 \\
        &  & 5.8266 & 80.703 & 23.334 \\ 
        \cline{3-5}
        &  & 0.00274 & 32.45 & 199.889 \\ 
        (C) $MSE_{R} < MSE_{P} < MSE_{OLS}$ & 469 (53.48\%) & 0.57782 & 76.68 & 3260.171 \\
        &  & 0.05795 & 99.5 & 4678.892 \\ 
        \cline{3-5}
        (D) $MSE_{R} < MSE_{OLS} < MSE_{P}$ & 0 (0\%) &  &  &  \\
        (E) $MSE_{OLS} < MSE_{P} < MSE_{R}$ & 0 (0\%) &  &  &  \\
        (F) $MSE_{OLS} < MSE_{R} < MSE_{P}$ & 0 (0\%) &  &  &  \\
        \hline
    \end{tabular}
    \caption{Comparison of mean squared errors for minimum coefficient of variation, maximum variance inflation factor and number of condition. First and third quartiles as well as the mean are provided.} \label{comparacionECM}
\end{table}

Table \ref{comparacionECM} compares the values obtained. It is observed that:
\begin{itemize}
    \item A total of 877 simulations were finally worked with. This is because there are 563 cases where it is not possible to calculate the minimum MSE for the ridge and/or penalized estimator either because the discretization is not fine enough or because the interval stop, 1, is insufficient. In this regard, tests have been performed by modifying the interval stop:
        \begin{itemize}
            \item If the interval stop is equal to 2: there are 903 cases in which both minima are obtained, in 432 cases (B) is verified and in 471 cases (C) is verified.
            \item IIf the interval stop is equal to 5: there are 921 cases in which both minima are obtained, in 445 cases (B) is verified and in 476 cases (C) is verified.
            \item If the interval stop is equal to 10: there are 923 cases in which both minima are obtained, in 447 cases (B) is verified and in 476 cases (C) is verified.
        \end{itemize}
        It can be seen that this increase results in a decreasing incorporation of simulations and that the new incorporations are mostly in case (B).
    \item Case B: In 46.52\% of the cases the penalized estimator has a lower mean square error and the ridge the worst. These cases are characterized by a high coefficient of variation, so that the multicollinearity treated is of the essential type.
    \item Case C: In 53.48\% of the cases the penalized estimator has a worse mean square error than the ridge, but better than that of OLS. These cases are characterized by a low coefficient of variation, so that the multicollinearity treated is of a non-essential type.
    \item The values shown by the condition number explain that the shift to the right discussed above is due to the detection by the condition number of non-essential multicollinearity, which is ignored by the variance inflation factor.
\end{itemize}

In short, the simulation shows that the penalized estimator performs better in terms of mean square error than the ridge estimator when the existing multicollinearity is of the essential type.

\begin{algorithm}
    \begin{algorithmic}[1]
        \REQUIRE Calculate $D(n)$ := \{discretization of the interval [0,1] with $n$ points\}
        \STATE i = 0
        \FOR {k $\in$ D(n)}
            \STATE i = i + 1
            \STATE Calculate $MSE \left( \widehat{\boldsymbol{\beta}}(k,h) \right)[i]$ with the expression given in subsection \ref{ecm_kh}
        \ENDFOR
        \STATE l = 0
        \FOR {k $\in$ 2:(i-1)}
            \IF {$MSE \left( \widehat{\boldsymbol{\beta}}(k,h) \right)[j] < MSE \left( \widehat{\boldsymbol{\beta}}(k,h) \right)[j-1]$ and $MSE \left( \widehat{\boldsymbol{\beta}}(k,h) \right)[j] < MSE \left( \widehat{\boldsymbol{\beta}}(k,h) \right)[j+1]$}
                \STATE l = l + 1
                \STATE index[l] = j
            \ENDIF
        \ENDFOR
        \STATE only = 0
        \IF {l = 1}
            \STATE only = 1
        \ENDIF
        \IF {l $>$ 1}
            \STATE only = 2
        \ENDIF
        \STATE MSE min = $MSE \left( \widehat{\boldsymbol{\beta}}(k,h) \right)[min(index)]$
    \end{algorithmic}
    \caption{Calculation of $k$ such that $MSE \left( \widehat{\boldsymbol{\beta}}(k,h) \right)$ is minimum} \label{algoritmo}
\end{algorithm}

\section{Augmented model and multicollinearity}
    \label{model.aumentado.SCRP}

From expression (\ref{est.SCRP.1}) the penalized estimator can be obtained by OLS estimation of the augmented model $\mathbf{y}_{A} = \mathbf{X}_{A} \cdot \boldsymbol{\beta} + \mathbf{u}_{A}$, where:
    $$\mathbf{X}_{A} = \left(
        \begin{array}{c}
            \mathbf{X} \\
            \sqrt{k} \cdot \mathbf{I}
        \end{array} \right)_{(n+p) \times p}, \quad
    \mathbf{y}_{A} = \left(
        \begin{array}{c}
            \mathbf{y} \\
            \sqrt{k} \cdot h \cdot \boldsymbol{\alpha}
        \end{array} \right)_{(n+p) \times 1}.$$
    However, the variance-covariance matrices would be different since:
    $$var \left( \widehat{\boldsymbol{\beta}}_{A} \right) = \sigma^{2} \cdot \left( \mathbf{X}^{t} \mathbf{X} + k \cdot \mathbf{I} \right)^{-1} \not= var \left( \widehat{\boldsymbol{\beta}}(k, h) \right).$$
    Similarly, the goodness-of-fit would also change since the dependent variable is different, $\mathbf{y} \not= \mathbf{y}_{A}$.

    It should be noted that the design matrix of the augmented model, $\mathbf{X}_{A}$, coincides with that of the augmented model of the ridge estimator (see Marquardt \cite{Marquardt1970}), then all operations based on them will coincide. Consequently, papers presented by García et al. \cite{Gracia2015} and Salmerón et al. \cite{Salmeron2017} to calculate, respectively, the extension of the Variance Inflation Factor (VIF) and the Condition Number (CN) to the ridge estimator are applicable to the penalized estimator. Thus, more specifically:
    \begin{itemize}
        \item García et al. \cite{Gracia2015} states that to obtain an extension of the VIF, $VIF(i,k)$ with $i=2,\dots,p$, with the desirable properties of continuity (for $k=0$ the calculated VIF coincides with that obtained in OLS), decreasing monotonicity as $k$ increases (since it is assumed that when applying the ridge estimate the degree of approximate multicollinearity decreases) and that a value greater than 1 is always obtained (since this is the minimum value of the VIF); data must be standardized\footnote{
                To standardize a set of data, center the data (subtract its mean) and divide it by its standard deviation multiplied by the square root of the number of observations.
            } and subsequently the model is augmented. In other words, to calculate the VIF for the ridge estimator, the following matrix must be used:
            $$\mathbf{x}_{A} = \left(
            \begin{array}{c}
                \mathbf{x} \\
                \sqrt{k} \cdot \mathbf{I}
            \end{array} \right)_{(n+p) \times p},$$
            where the matrix $\mathbf{x}$ is the matrix resulting from the standardization of the columns of $\mathbf{X}$.
        \item Salmerón et al. \cite{Salmeron2017} it is recommended that to obtain an extension of the CN, $CN(k)$, that verifies the above desirable properties (continuity at $k=0$, decreasing monotonicity at $k$ and being greater than 1 for any value of $k$) it has to be calculated from the expression:
            $$CN(k) = \sqrt{\frac{\xi_{max} + k}{\xi_{min} + k}},$$
             where$\xi_{min}$ and $\xi_{max}$ are, respectively, the minimum and maximum eigenvalues of $\mathbf{X}_{T}^{t} \mathbf{X}_{T}$ where $\mathbf{X}_{T}$ denotes the matrix $\mathbf{X}$ transformed. These authors recommend that the transformation should be typing\footnote{
                To typify a set of data it is necessary to center it (subtract its mean) and divide it by its standard deviation.
            } the data, since in such a case the results obtained are similar to those obtained from the VIF. However, this transformation ignores nonessential type multicollinearity (see, for example, Salmerón et al. \cite{salmeron2018JSCS,Salmeron2020}), so this transformation is not going to be used in this paper. In accordance with the recommendations of Belsley et al. \cite{Belsleyetal1980} or Belsley \cite{Belsley1984,Belsley1982} the eigenvalues of the matrix will be calculated. $\mathbf{X}_{T}^{t} \mathbf{X}_{T}$ where the transformation used is that of unit length\footnote{
                To transform a set of data into unit length, divide it by the square root of the sum of each element squared. Note that if the data have zero mean (a property that is obtained by centering them), it is verified that this transformation coincides with the standardization.
            }.
    \end{itemize}

\section{The penalty parameter $k$} \label{elec.SCRP}

In order to obtain the value of the above expressions it is necessary to use a concrete value of the penalty parameter, $k$. For this reason, some possible criteria for determining a concrete value of the parameter, as well as its interpretation, are proposed below.

    \subsection{Interpretation of the penalty parameter}

    Taking into account Remark \ref{kas} in which $k = \frac{k_{2}}{k_{1}}$ for $k_{1} > 0$ and considering that $k_{1} + k_{2} = 1$, it is possible to obtain the relative importance of each of the two summands of the expression (\ref{lagran.3}) in the minimization performed.

    That is, assuming that $k = r \geq 0$ it is obtained thaat $k_{1} = \frac{1}{1+r}$ and $k_{2} = \frac{r}{1+r}$. Thus, for example:
    \begin{itemize}
        \item If $k = 1$, then $k_{1} = 0.5$ and $k_{2} = 0.5$. In this case, the same importance would be given to $SCR(\boldsymbol{\beta})$ and to $(\boldsymbol{\beta} - h \cdot \boldsymbol{\alpha})^{t} (\boldsymbol{\beta} - h \cdot \boldsymbol{\alpha})$.
        \item If $k = 2$, then $k_{1} = \frac{1}{3}$ and $k_{2} = \frac{2}{3}$. In this case, double importance is given to $(\boldsymbol{\beta} - h \cdot \boldsymbol{\alpha})^{t} (\boldsymbol{\beta} - h \cdot \boldsymbol{\alpha})$ que a $SCR(\boldsymbol{\beta})$.
    \end{itemize}

    In general:
    \begin{itemize}
        \item If $k < 1$ it is obtained that $k_{2} < k_{1}$ and, then, it would be giving more relevance to $SCR(\boldsymbol{\beta})$ than$(\boldsymbol{\beta} - h \cdot \boldsymbol{\alpha})^{t} (\boldsymbol{\beta} - h \cdot \boldsymbol{\alpha})$.
        \item If $k > 1$ it is obtained that $k_{2} > k_{1}$ and, then, it would be giving more relevance to a $(\boldsymbol{\beta} - h \cdot \boldsymbol{\alpha})^{t} (\boldsymbol{\beta} - h \cdot \boldsymbol{\alpha})$ than $SCR(\boldsymbol{\beta})$.
    \end{itemize}
    Note that in the ridge estimator the value of $k$ has traditionally been considered to belong to the interval $[0, \ 1]$. hat is, preference is given to minimizing the sum of squares of the residuals over the shrinkage of the coefficients.

    Moreover, the fact that $k_{1}$ is imposed to be greater than zero implies that the SSR is always minimized. At the same time, since the higher the value of $k$ the lesser the role of $SCR(\boldsymbol{\beta})$, it is intuited that the goodness-of-fit calculated from $GoF(k, h)$ will be decreasing..

    \subsection{Selection of the penalty parameter}
        \label{elegir_k}

    A first option, following  Hoerl and Kennard \cite{HKa, HKb}, is to select that value of $k$ that stabilizes the estimates of the coefficients of the independent variables of the model. For this purpose, it is necessary to plot the values of  $\widehat{\boldsymbol{\beta}}(k, h)$ as a function of $k$ nd observe at what value of $k$ the estimates are stable.

    On the other hand, following García et al \cite{Gracia2015} and Salmerón et al \cite{Salmeron2017}, could be considered the first value of $k$ that makes the extension of VIF or CN below the thresholds (10 and 20, respectively) above which multicollinearity is traditionally considered to be troubling.

    Furthermore, since the goal is to make the estimates similar to those obtained in $\boldsymbol{\alpha}$, could be considered the value of $k$ that makes $\widehat{\boldsymbol{\beta}}(k, h)$ and $\boldsymbol{\alpha}$ close to each other. For example, the difference between the two could be imposed to be less than 10\%:
    $$\frac{|| \boldsymbol{\alpha} - \widehat{\boldsymbol{\beta}}(k, h) ||}{|| \boldsymbol{\alpha} ||} < 0.1.$$
    Section \ref{rate} of Appendix \ref{apen.operaciones} shows that $|| \boldsymbol{\alpha} - \widehat{\boldsymbol{\beta}}(k, h) ||$ decreases toward zero when $k$ increases and $h=1$.

    Finally, by following Hoerl et al \cite{HKB1975}, the value of $k$ that minimizes the mean square error (if it exists) could be considered..

    \begin{example}[continuation of Example \ref{ejemplo.SCRP}]
        \label{ejemplo.SCRP.2}
        In example \ref{ejemplo.SCRP} (where the degree of essential multicollinearity is troubling) considering  that $k \in \{0, 0.01, 0.02, \dots, 99.99, 100 \}$ and $h=1$, the following values of $k$ are obtained by considering different criteria:
        \begin{itemize}
            \item $k=0.06$ is the value of $k$ for which the maximum $VIF(i,k)$, $i=2,3$, is less than 10 (concretely 9.060921);
            \item $k=0.01$ is the value of $k$ for which the $CN(k)$ is lower that 20 (concretely 14.2435) and $k=0.01$ is the value of $k$ for which the $CN(k)$ lower that 10 (concretely 8.316654)
            \item $k=0.03$ is the value of $k$ for which the mean square error is minimum, concretely is equal to 4.284018 taht is lower that the one obtained by OLS 4.311986.
        \end{itemize}
        On the other hand, in the considered discretization of the interval [0, 100] there does not exist $k$ such that  $\frac{|| \boldsymbol{\alpha} - \widehat{\boldsymbol{\beta}}(k, h) ||}{|| \boldsymbol{\alpha} ||} < 0.1$. The minimum value obtained in this case is 0.4402975.

        Finally, with these values of $k$ we would be giving more weight to the minimization of $SCR(\boldsymbol{\beta})$ than to $(\boldsymbol{\beta} - \boldsymbol{\alpha})^{t} (\boldsymbol{\beta} - \boldsymbol{\alpha})$ according to the values shown in the following table:
        \begin{center}
            \begin{tabular}{ccc}
                \hline
                $k$ & $k_{1}$ & $k_{2}$ \\
                \hline
                0.01 & 0.990099  & 0.00990099  \\
                0.03 & 0.9708738  &  0.02912621 \\
                0.06 & 0.9433962  &  0.05660377 \\
                \hline
            \end{tabular}
        \end{center}
        \hfill $\Box$
    \end{example}


\section{Inference}
    \label{inf.bootstrap}

Although there are various attempts to address inference in the ridge estimator (see, e.g., Obenchain \cite{Obenchain1975,Obenchain1977}, Halawa and El Bassiouni \cite{HalawaBassiouni200},  Gökpınar and Ebegil \cite{GokpinarEbegil2016}, Sengupta and Sowell \cite{SenguptaSowell2020} or Vanhove \cite{Vanhove2021}), in this paper we will focus on the use of bootstrap methods (see, for example, Efron \cite{Efron1979,Efron1981,Efron1982,Efron1987}, Efron and Gong \cite{EfronGong1983} or the a review of the above-mentioned works presented by Efron and Tibshirani \cite{EfronTibshirani1986}).

Thus, given a fixed value of $k$ obtained by one of the methods proposed in subsection \ref{elegir_k} and a value of $\boldsymbol{\alpha}$, the following steps will be taken:
\begin{enumerate}[(i)]
    \item Generate randomly and with replacement $m$ subsamples of equal size to the original one. The value of $m$ must be large.
    \item For each sub-sample above, the statistic of interest $theta$ is calculated. Therefore, we have $m$ values for this statistic: $\theta_{1},\dots,\theta_{m}$.
    \item Calculate the mean, $\overline{\theta} = \frac{1}{m} \sum \limits_{i=1}^{m} \theta_{i}$, and sample quasi-variance, $\sigma^{2}_{\theta} = \frac{1}{m-1} \sum \limits_{i=1}^{m} (\theta_{i} - \overline{\theta})^{2}$, for the $m$ values obtained in previous step. Both measures can be used to obtain the approximation of a confidence interval by using expression $\overline{\theta} \pm 1.96 \cdot \sigma_\theta$. \\
        Similarly, $\overline{\theta}$ can be used as an approximation of the point estimator of  $\theta$.
    \item Another confidence region can be obtained alternatively by considering the lower and upper percentiles, 0.025 and 0.975, as the lower and upper extremes respectively, of the $m$ values calculated in the second step: $[P_{0.025}(\theta_{1},\dots,\theta_{m}), P_{0.975}(\theta_{1},\dots,\theta_{m})]$.
\end{enumerate}
EWhat is interesting in this work are the cases in which $\theta$ is equal to $\widehat{\boldsymbol{\beta}}(k,1)$ or $GoF(k,1)$.

\begin{example}[continuation of the Example \ref{ejemplo.SCRP.2}]
    \label{ejemplo.SCRP.3}

   Considering:
    \begin{itemize}
        \item values 0.01, 0.03 and 0.06 of $k$ which make, respectively, that $CN(k)<20$; the mean square error is minimum and $CN(k)<10$; or $VIF(i,k)<10$ for $i=2,3$ (see Example \ref{ejemplo.SCRP.2}),
        \item $\boldsymbol{\alpha} = (-29.476251, -18.036792, -3.601057)^{t}$ y
        \item $m=10000$,
    \end{itemize}
    Table \ref{tab:bootstrap} shows the approximation of the point and interval estimators of $\widehat{\boldsymbol{\beta}}(k,1)$ and $GoF(k,1)$.

    It can be seen that the approximation of the point estimators, $\overline{\theta}$, are very close to the estimators of $\widehat{\boldsymbol{\beta}}(k,1)$ and calculation of $GoF(k,1)$ obtained from expressions (\ref{est.SCRP.1}) and (\ref{ba.ajuste.alterna}), respectively. This facilitates that the values of  $\widehat{\boldsymbol{\beta}}(k,1)$ and $GoF(k,1)$ belong to the two calculated confidence intervals.

    Moreover, from these ranges it could be argued that the coefficient estimates $\beta_{1}$ and $\beta_{3}$ are significantly different from zero for the three values of $k$ considered since the intervals obtained do not contain zero, while it cannot be discarded that the estimate of the $k$ coefficient $\beta_{2}$ is null.
    The same behaviour is observed as in OLS, perhaps due to the fact that the estimates obtained are still far from the values for $\boldsymbol{\alpha}$.

    In this context, considering that $k=100$ (maximum value of $k$ considered in Example \ref{ejemplo.SCRP.2} to determine a value for this parameter), it is observed that the estimates of $\widehat{\boldsymbol{\beta}}(k,1)$ change substantially from the previous ones. In this case, all the coefficients can be considered significantly different from zero, although the estimation of the constant and the third coefficient are far from the values proposed in $\boldsymbol{\alpha}$.

    It is curious to note that there has been a distancing of $\widehat{\boldsymbol{\beta}}_{3}(k,1)$ from -3.601057: the value -3.6767991 for $k=0.01$ changes to $0.3170817$ for $k=100$. Note hat for $k=100$ a weight of 0.00990099 would be given to minimisation of $SCR(\boldsymbol{\beta})$ against a value of 0.990099 for  $(\boldsymbol{\beta} - \boldsymbol{\alpha})^{t} (\boldsymbol{\beta} - \boldsymbol{\alpha})$.

    Finally, it is considered that $k=1000$, the estimation of $\widehat{\boldsymbol{\beta}}_{3}(k,1)$ is equal to 0.03183904 with approximations of the associated confidence intervals equal to $[-0.06174285,   0.1239955]$ and $[-0.0678132,   0.118499]$, i.e., a coefficient not significantly different from zero is obtained. Once again, it is interesting to note that the third variable initially has values for  $\widehat{\beta}_{3}$ and $\widehat{\alpha}_{3}$ very similar between them (see Example \ref{ejemplo.SCRP}) and that as a consequence of the linear relationship with the second variable, the coefficient becomes not significantly different from zero using the penalized estimator.
    \hfill $\Box$
\end{example}

\begin{sidewaystable}
    \centering
    \begin{tabular}{cccccc}
        \hline
        $i$ & $\theta = \widehat{\boldsymbol{\beta}}_{i}(k,1)$, $k=0.01$ & $\overline{\theta}$ & $\sigma_\theta$ & $\overline{\theta} \pm 1.96 \cdot \sigma_\theta$ & $[P_{0.025}(\theta_{1},\dots,\theta_{m}), P_{0.975}(\theta_{1},\dots,\theta_{m})]$ \\
        \hline
        1 & 4.3580904 & 4.3637278 & 0.2876449 & [3.803312,  4.924143] & [3.790357,  4.914893] \\
        2 & 0.3809489 & 0.3793928 & 0.9808815 & [-1.550971,  2.309756] & [-1.564748,  2.283066] \\
        3 & -3.6770687 & -3.6767991 & 0.1959585 & [-4.062437, -3.291161] & [-4.057372, -3.287766] \\
        \hline
        $i$ & $\theta = \widehat{\boldsymbol{\beta}}_{i}(k,1)$, $k=0.03$ & $\overline{\theta}$ & $\sigma_\theta$ & $\overline{\theta} \pm 1.96 \cdot \sigma_\theta$ & $[P_{0.025}(\theta_{1},\dots,\theta_{m}), P_{0.975}(\theta_{1},\dots,\theta_{m})]$ \\
        \hline
        1 & 4.3264551 & 4.3226154 & 0.281187 & [3.757987,  4.887244] & [3.743654,  4.876066] \\
        2 & 0.2552237 & 0.2439003 & 0.9683303 & [-1.667152,  2.154952] & [-1.708548,  2.117196] \\
        3 & -3.6519522 & -3.6496978 & 0.1934983 & [-4.031710, -3.267686] & [-4.025403, -3.259135] \\
        \hline
        $i$ & $\theta = \widehat{\boldsymbol{\beta}}_{i}(k,1)$, $k=0.06$ & $\overline{\theta}$ & $\sigma_\theta$ & $\overline{\theta} \pm 1.96 \cdot \sigma_\theta$ & $[P_{0.025}(\theta_{1},\dots,\theta_{m}), P_{0.975}(\theta_{1},\dots,\theta_{m})]$ \\
        \hline
        1 & 4.2795271 & 4.27471093 & 0.2873334 & [3.717146,  4.832276] & [3.695980,  4.823938] \\
        2 & 0.0691914 & 0.05400399 & 0.9622777 & [-1.831101,  1.939109] & [-1.810631,  1.918534] \\
        3 & -3.6147874 & -3.61174146 & 0.1923354 & [-3.988318, -3.235165] & [-3.984625, -3.238844] \\
        \hline
        $i$ & $\theta = \widehat{\boldsymbol{\beta}}_{i}(k,1)$, $k=100$ & $\overline{\theta}$ & $\sigma_\theta$ & $\overline{\theta} \pm 1.96 \cdot \sigma_\theta$ & $[P_{0.025}(\theta_{1},\dots,\theta_{m}), P_{0.975}(\theta_{1},\dots,\theta_{m})]$ \\
        \hline
        1 & -14.7479306 & -14.8279630 & 0.3758633 & [-15.564655, -14.0912709] & [-15.6807864, -14.1957101] \\
        2 & -19.3578449 & -19.3543196 & 0.2441033 & [-19.832762, -18.8758772] & [-19.8351704, -18.8815710] \\
        3 & 0.3169069 & 0.3170817 & 0.05212487 & [0.214917,   0.4192465] & [0.2151219,   0.4176053] \\
        \hline
        $k$ & $\theta = GoF(k,1)$ & $\overline{\theta}$ & $\sigma_\theta$ & $\overline{\theta} \pm 1.96 \cdot \sigma_\theta$ & $[P_{0.025}(\theta_{1},\dots,\theta_{m}), P_{0.975}(\theta_{1},\dots,\theta_{m})]$ \\
        \hline
        0.01 & 0.999884 & 0.9998847 & 2.468503$\cdot 10^{-5}$ & [0.9998363, 0.9999331] & [0.9998286, 0.9999254] \\
        0.03 & 0.999884 & 0.9998847 & 2.476162$\cdot 10^{-5}$ & [0.9998362, 0.9999333] & [0.9998297, 0.9999260] \\
        0.06 & 0.9998838 & 0.9998844 & 2.470472$\cdot 10^{-5}$ & [0.9998360, 0.9999329] & [0.9998289, 0.9999251] \\
        100 & 0.992539 & 0.9924073 & 0.001157019 & [0.9901396, 0.9946751] & [0.9898363, 0.9943386] \\
        \hline
    \end{tabular}
    \caption{Bootstrap inference for the second case in the Example \ref{ejemplo.SCRP}} \label{tab:bootstrap}
\end{sidewaystable}

\section{Stability of the penalized estimator}
    \label{stability}

    In order to study the stability of the estimators, the independent variables of the multiple linear regression model will be perturbed by 1\% following the procedure shown in Belsley \cite{Belsley1984}. Thus, given a vector $\mathbf{x}$ of size $n$, a perturbation of it by 1\% is given by the expression::
    $$\mathbf{x}_{p} = \mathbf{x} + 0.01 \cdot \mathbf{p} \cdot \frac{||\mathbf{x}||}{||\mathbf{p}||},$$
    where $\mathbf{p}$ is a random vector, of equal dimension to $\mathbf{x}$, and $||\mathbf{x}|| = \sqrt{\sum \limits_{i=1}^{n} x_{i}^{2}}$.

    The effect of the perturbation on the estimation of the coefficients is quantified by the following expression:
    \begin{equation}
        \frac{||\widetilde{\boldsymbol{\beta}} - \widetilde{\boldsymbol{\beta}}_{p}||}{||\widetilde{\boldsymbol{\beta}}||} \cdot 100\%,
        \label{cuantification_perturb}
    \end{equation}
    where $\widetilde{\boldsymbol{\beta}}$ denotes an estimator of the parameter $\boldsymbol{\beta}$ and the subindex $p$ refers to the estimates obtained in the perturbed model, i.e. in the model in which all the independent variables have been perturbed.\footnote{
        In this case it is considered that all variables are quantitative. If there were any binary variables, it would be advisable to use a different perturbation method than the one described.
    } (without considering the intercept).

    To obtain consolidated results, the above procedure will be performed 1000 times and its mean and confidence region determined by the 0.025 and 0.975 percentiles will be used for interpretations.

\begin{figure}
  \centering
  \includegraphics[width=8cm]{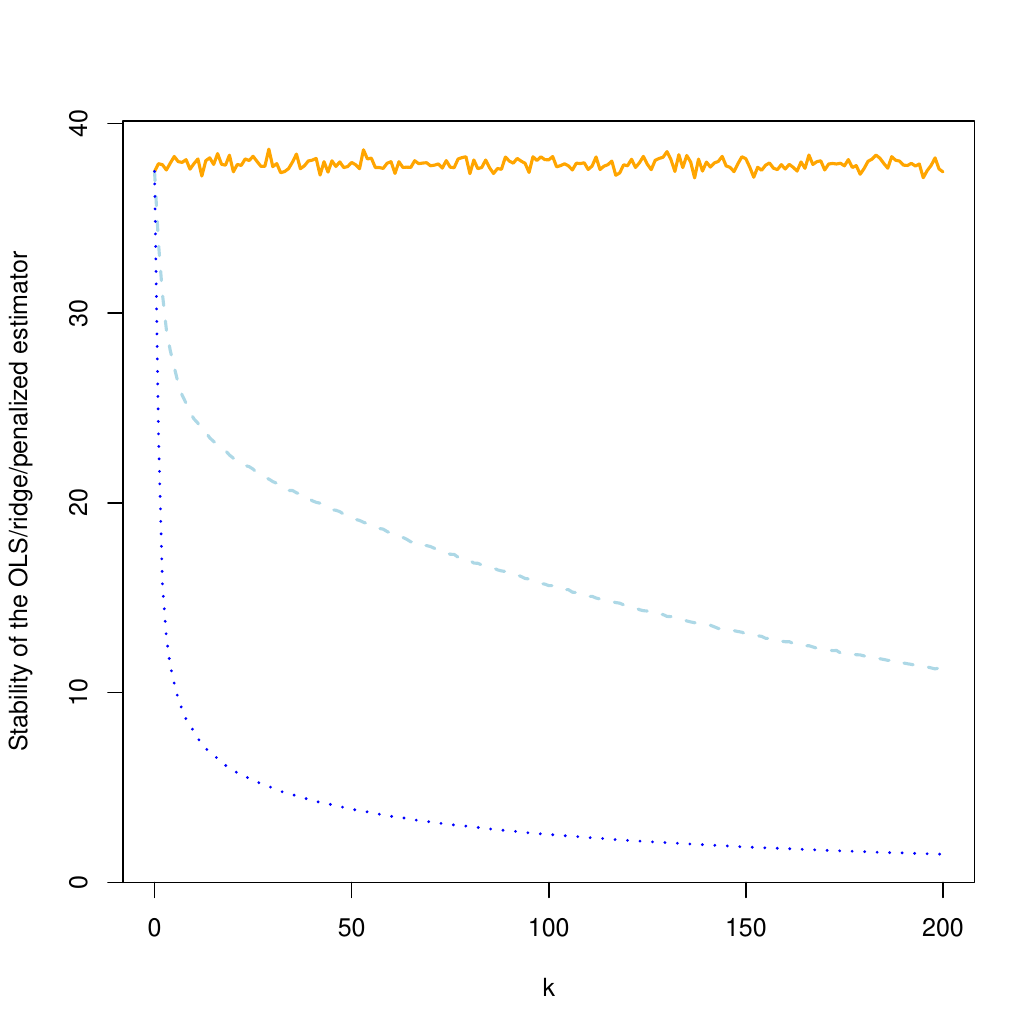}
  \caption{Average percentage change (for 1000 iterations) that occurs in the coefficient estimates in the presence of a 1\% perturbation in the independent variables for  $k \in \{ 0, 1, 2, \dots, 198, 199, 199 \}$: OLS estimator in orange, ridge estimator in lightblue and penalized estimator in blue}\label{fig:norms_example}
\end{figure}

\begin{example}[continuation of example \ref{ejemplo.SCRP.3}]
    \label{ejemplo.SCRP.4}

   In example \ref{ejemplo.SCRP} in which the degree of essential approximate multicollinearity is troubling, the independent variables have been perturbed 1000 times for $k=0$ and for each value of $k$ considered in the Example \ref{ejemplo.SCRP.3}: $k = 0.01, 0.03, 0.06, 100, 1000$.

    For $k=0$ (i.e., when OLS is applied), is obtained that perturbing the independent variables by 1\% implies a variation of 38.27699\%. Therefore, this would be a case in which the multicollinearity in the model causes numerical instability, since a small change in the independent variables implies a significant variation in the estimates of their coefficients.

     Table \ref{tab:estabilidad} shows the mean, 0.025 percentile and 0.095 percentile of the percent variance of the 1000 iterations performed for $k = 0, 0.01, 0.03, 0.06, 100, 1000$ and for each estimator: OLS, ridge and penalized:
    \begin{itemize}
        \item When $k = 0.01, 0.03, 0.06$, it is observed that the values obtained for both the mean value and the percentiles are very similar to those of OLS, i.e., those values of $k$ that make the quadratic error minimal or measures such as the condition number or variance inflation factor fall below the established thresholds do not substantially improve the numerical instability of the model.
        \item It is observed that for the high values of $k$ considered, $k = 100, 1000$, stable estimators are achieved in the presence of small changes in the independent variables.
    \end{itemize}

    To finish, Figure \ref{fig:norms_example} shows the average percentage change (out of 1000 iterations) that occurs in the coefficient estimates in the presence of a 1 percent perturbation in the independent variables for   $k \in \{ 0, 1, 2, \dots, 198, 199, 199 \}$. It is observed that the penalized estimator (blue - dotted) has a higher stability to the discussed perturbations than the ridge estimator (lightblue - dashed) as the value of $k$ increases.
    \hfill $\Box$
\end{example}

\begin{table}
    \centering
    \begin{tabular}{cccccccccc}
        \hline
         &  & Ridge &   &  & Penalized &  \\
        $k$ & Mean & $P_{0.025}$ & $P_{0.075}$  & Mean & $P_{0.025}$ & $P_{0.075}$ \\
        \hline
        0 & 38.27699 & 23.95753 & 55.43604  & 38.27699 & 23.95753 & 55.43604 \\
        0.01 & 38.22262 & 23.95097 &  55.33613 & 38.20928 & 24.02455 & 55.36458 \\
        0.03 & 38.11476 & 23.94311 & 55.13764   & 38.06461 & 24.15165 & 55.20508 \\
        0.06 & 37.95507 & 23.93155 & 54.84321  & 37.82431 & 24.26060 & 54.92455 \\
        100 & 15.65896 & 13.62525 & 17.41957  & 2.535079 & 2.156762 & 2.913410 \\
        1000 & 3.256717 & 2.822840 & 3.643900  & 0.3383532 & 0.3016403 & 0.3746771 \\
        \hline
    \end{tabular}
    \caption{Mean, 0.025th percentile and 0.095th percentile of the percent variation of the 1000 iterations performed for $k = 0, 0.01, 0.03, 0.06, 100$ in each estimator.} \label{tab:estabilidad}
\end{table}

\section{Example} \label{ejemplos.SCRP}

    \begin{table}
        \begin{center}
        \begin{tabular}{ccccccc}
            \hline
            Year & $\mathbf{D}$ & $\mathbf{C}$ & $\mathbf{I}$ \\
            \hline
            1996 & 3.80510 & 4.7703 & 4.8786  \\
            1997 & 3.94580 & 4.7784 & 5.0510  \\
            1998 & 4.05790 & 4.9348 & 5.3620  \\
            1999 & 4.19130 & 5.0998 & 5.5585  \\
            2000 & 4.35850 & 5.2907 & 5.8425  \\
            2001 & 4.54530 & 5.4335 & 6.1523  \\
            2002 & 4.81490 & 5.6194 & 6.5206  \\
            2003 & 5.12860 & 5.8318 & 6.9151  \\
            2004 & 5.61510 & 6.1258 & 7.4230  \\
            2005 & 6.22490 & 6.4386 & 7.8024  \\
            2006 & 6.78640 & 6.7394 & 8.4297  \\
            2007 & 7.49440 & 6.9104 & 8.7241  \\
            2008 & 8.39930 & 7.0993 & 8.8819  \\
            2009 & 9.39510 & 7.2953 & 9.1636  \\
            2010 & 10.68000 & 7.5614 & 9.7272  \\
            2011 & 12.07100 & 7.8036 & 10.3010  \\
            2012 & 13.44821 & 8.0441 & 10.9830  \\
            \hline
        \end{tabular}
        \caption{Data on credit in the United State} \label{datos.Wissel}
        \end{center}
    \end{table}

Table \ref{datos.Wissel} shows the outstanding mortgage debt, $\mathbf{D}$, personal consumption, $\mathbf{C}$, personal income, $\mathbf{I}$, and outstanding consumer credit, $\mathbf{CP}$, for years from 1996 to 2012 in relation to credit in the United State. This dataset was previously applied by \cite{Wissell}. Results of the OLS estimation are shown in Table \ref{tab.OLS_Wissel}.

\begin{table}
    \centering
    \begin{tabular}{ccc}
        \hline
        Variable & Estimate & p-value \\
        \hline
        \multirow{2}{*}{Intercept} & 5.469264 & \multirow{2}{*}{0.681} \\
         & (13.016791) &  \\
        \multirow{2}{*}{$\mathbf{C}$} & -4.252429 & \multirow{2}{*}{0.423} \\
         & (5.135058) &  \\
        \multirow{2}{*}{$\mathbf{I}$} & 3.120395 & \multirow{2}{*}{0.149} \\
         & (2.035671) &  \\
        \multirow{2}{*}{$\mathbf{CP}$} & 0.002879 & \multirow{2}{*}{0.626} \\
         & (0.005764) &  \\
        F-statistic & 52.3 & 0.0000001629 \\
        $R^{2}$ & \multicolumn{2}{c}{0.9235} \\
        \hline
    \end{tabular}
    \caption{OLS estimation of dataset related to bank credit in the United State} \label{tab.OLS_Wissel}
\end{table}

It is observed that the OLS estimates of the coefficients on consumption, income and credit outstanding are not significantly\footnote{Throughout the section, it is considered to work at 5\% significance} different from zero, while at the same time the model is found to be jointly valid. This contradiction suggests that the degree of approximate multicollinearity in the model affects the statistical analysis of the model.

For this model, the following simple regressions can be performed:
$$\mathbf{D} = \alpha_{1} + \alpha_{2} \mathbf{C} + \mathbf{v}, \quad
\mathbf{D} = \alpha_{1} + \alpha_{3} \mathbf{I} + \mathbf{v}, \quad
\mathbf{D} = \alpha_{1} + \alpha_{4} \mathbf{CP} + \mathbf{v},$$
which in all cases have coefficients for the slopes significantly different from zero and lead to establish $\boldsymbol{\alpha} = \left( 6.76245, 2.62875, 1.51533, 0.00519 \right)^{t}$.

\subsection{Diagnosis of troubling approximate multicollinearity.}

This possible worrying approximate multicollinearity is confirmed with the conclusions obtained from the following diagnosis measures (calculated by ``multiColl'' package \cite{multiColl} R \cite{RCoreTeam}):
\begin{itemize}
    \item the coefficients of variation of the independent variables are equal to 0.1718940, 0.2482804 and 0.3607848 (not lower than the threshold of  0.1002506 established as troubling by Salmerón et al. \cite{Salmeron2020}),
    \item the values for the variance inflation factor are equal to 589.7540, 281.8862 and 189.4874 (higher than the threshold of 10 established as troubling, for example, by Marquardt  \cite{Marquardt1970}),
    \item the condition number is equal to 332.3 (above the threshold of 30 established as troubling by Belsley  et al. \cite{Belsleyetal1980} or Belsley \cite{Belsley1984,Belsley1982}) and
    \item the determinant of the correlation matrix is equal to 0.00002007699 (below the threshold of 0.1013 established as troubling by Garcia et al. \cite{Garcia2018}).
\end{itemize}
More specifically, it can be stated that the degree of essential multicollinearity existing in model is troubling, while that of the non-essential type is not.

\subsection{Traces}

\begin{figure}
  \centering
  \includegraphics[width=8cm]{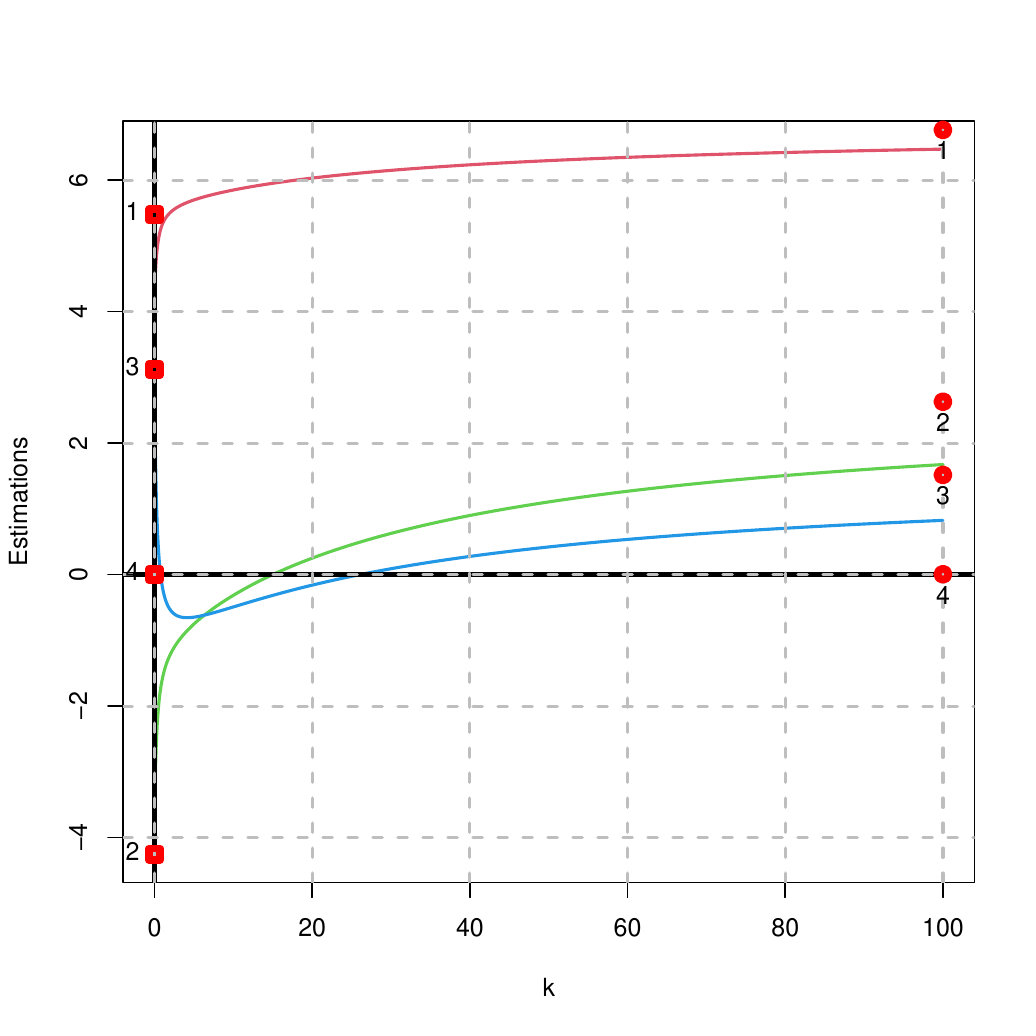}
  \includegraphics[width=8cm]{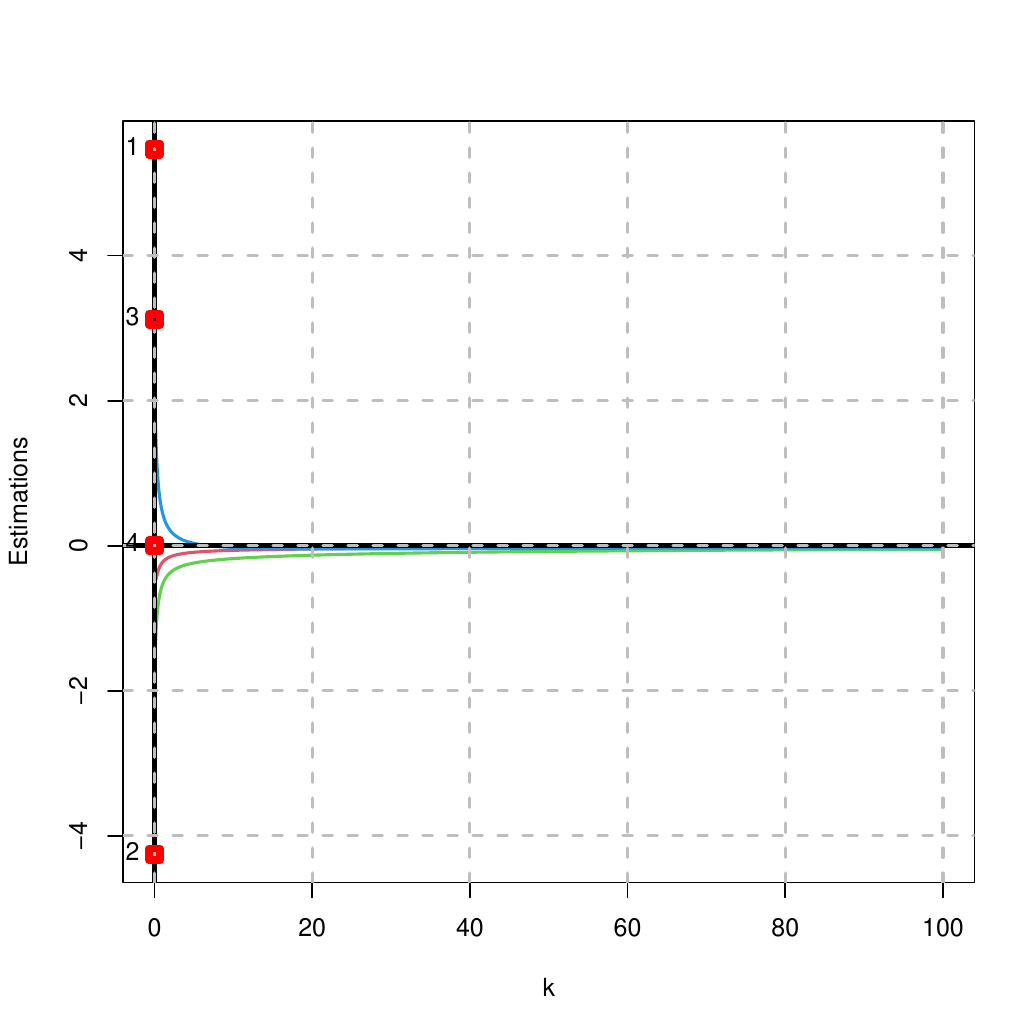}
  \caption{Trace of the coefficient estimates of the penalized estimator ($h=1$, left) and ridge ($h=0$, right) of the model on bank credit in the United States for $k \in \{ 0, 0.01, 0.02, \dots, 100 \}$} \label{fig.example.1}
\end{figure}

\begin{figure}
  \centering
  \includegraphics[width=8cm]{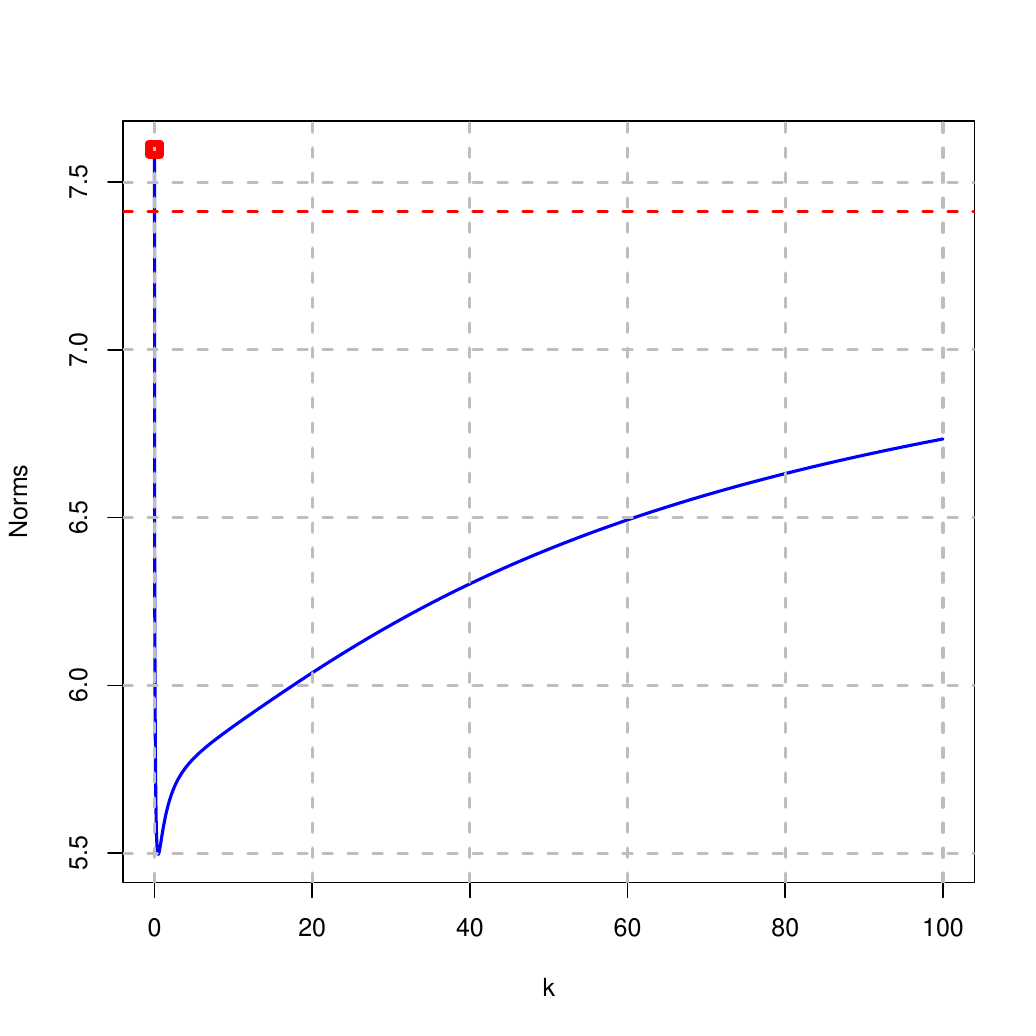}
  \includegraphics[width=8cm]{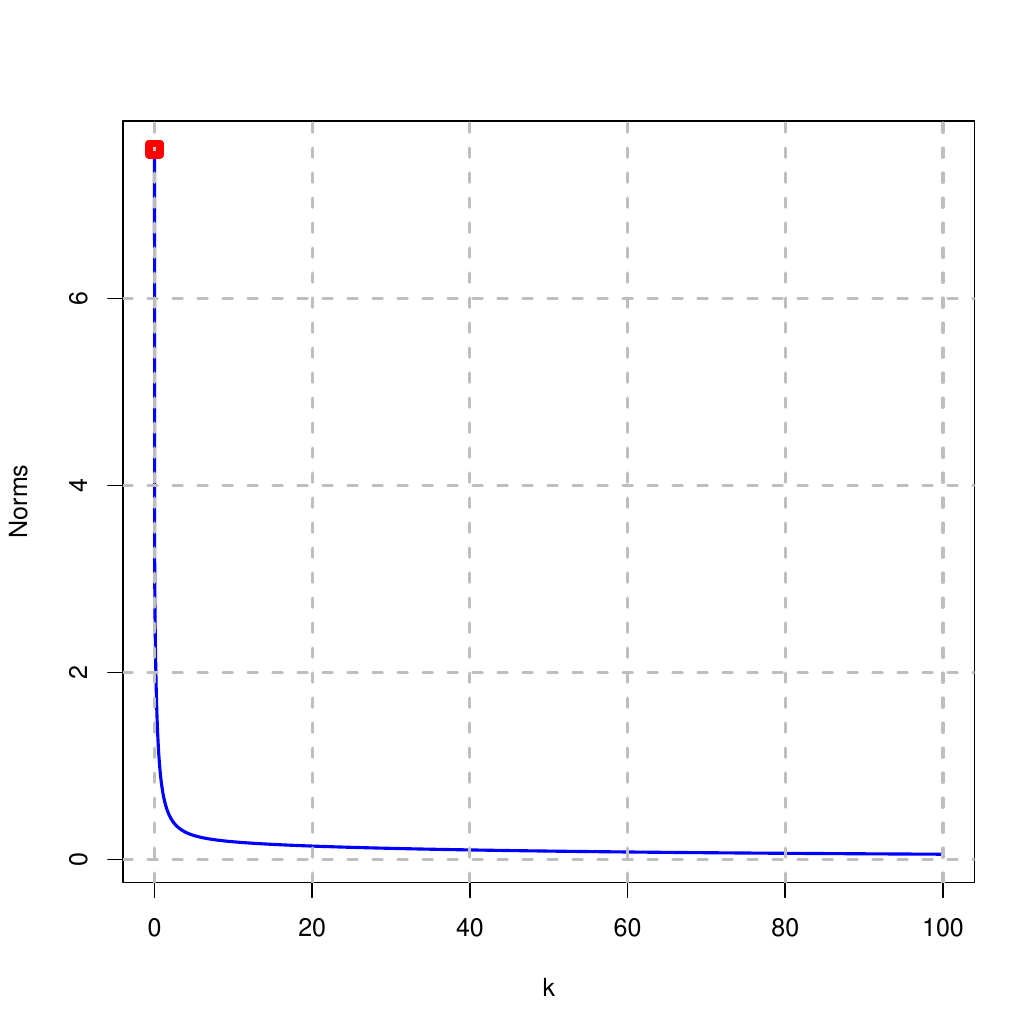}
  \caption{Trace of the norms of the coefficient estimates of the penalized estimator ($h=1$, left) and ridge ($h=0$, right) of the model on bank credit in the United States for $k \in \{ 0, 0.01, 0.02, \dots, 100 \}$} \label{fig.example.2}
\end{figure}

\begin{figure}
  \centering
  \includegraphics[width=8cm]{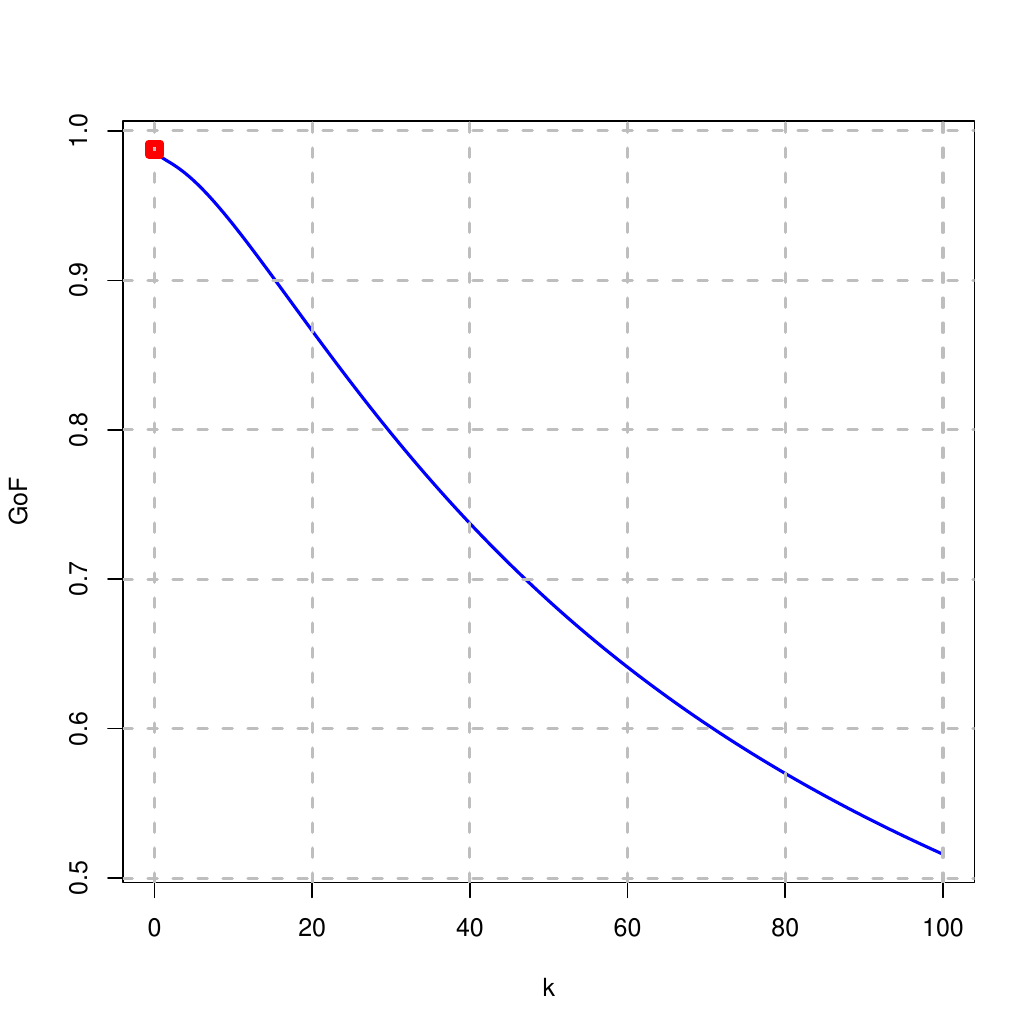}
  \includegraphics[width=8cm]{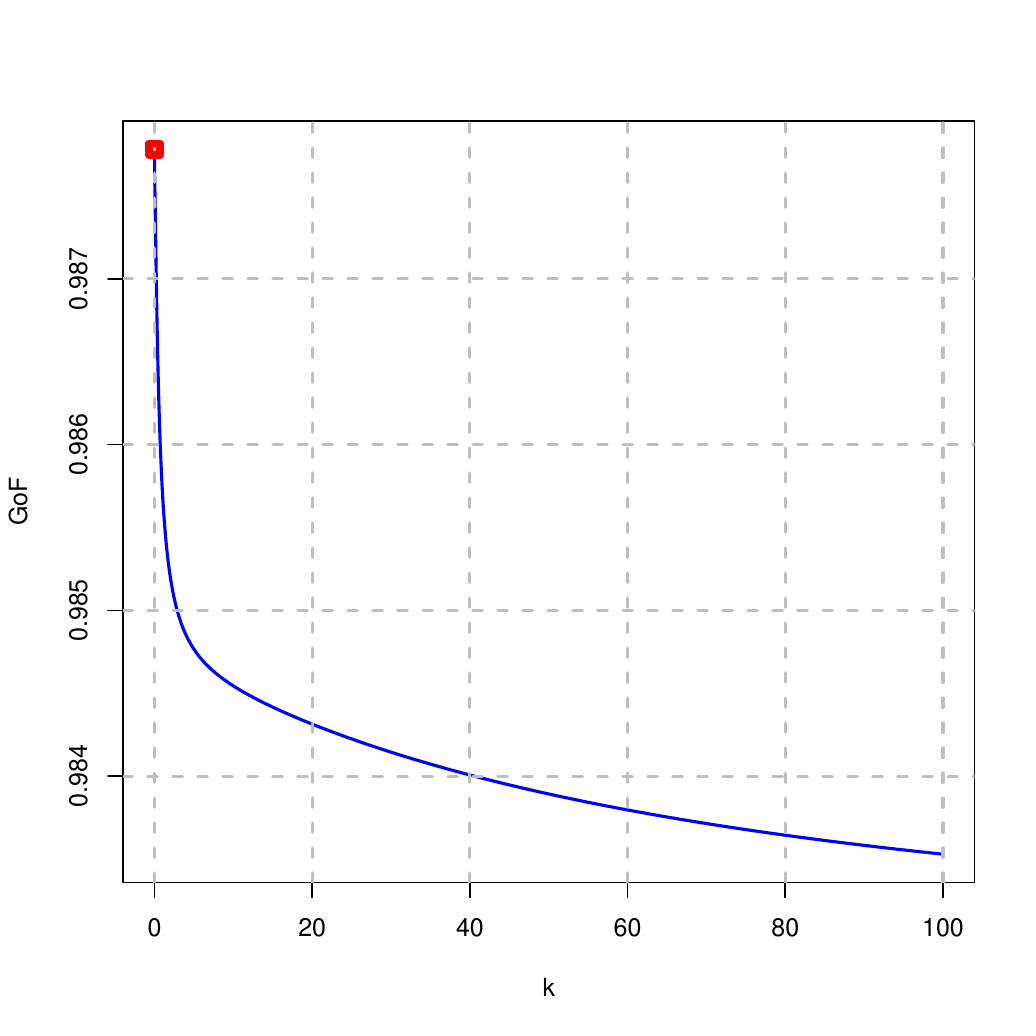}
  \caption{Trace of the goodness-of-fit of the penalized estimator ($h=1$, left) and ridge ($h=0$, right) of the model on bank credit in the United States for $k \in \{ 0, 0.01, 0.02, \dots, 100 \}$} \label{fig.example.3}
\end{figure}

Figures \ref{fig.example.1} to \ref{fig.example.3} show the traces of the estimates, mean squared error, goodness of fit and norms on the difference of the penalized estimator (right) and ridge (left) considering that $k \in \{ 0, 0.01, 0.02, \dots, 100 \}$. It is concluded that:
\begin{itemize}
    \item The estimates of the ridge estimator decay rapidly towards zero; while that of the penalized estimator do not, as they converge to the values given for the vector $\boldsymbol{\alpha}$.
    \item The norm of the penalized estimator decreases until it reaches a minimum value, from which it grows in the direction of its horizontal asymptote; while the norm of the ridge estimator is clearly strictly decreasing towards zero.
    \item The goodness of fit in both cases is decreasing, although in the case of the penalized estimator it takes smaller values than in the ridge (it decreases more rapidly).
    \item It is observed that in both cases there is a value of $k$ that implies a minimum squared error that is smaller than that of OLS. Moreover, it is found that the mean squared error of the penalized and ridge estimator is smaller than that of OLS for all positive values of $k$ considered.
\end{itemize}

\begin{figure}
  \centering
  \includegraphics[width=8cm]{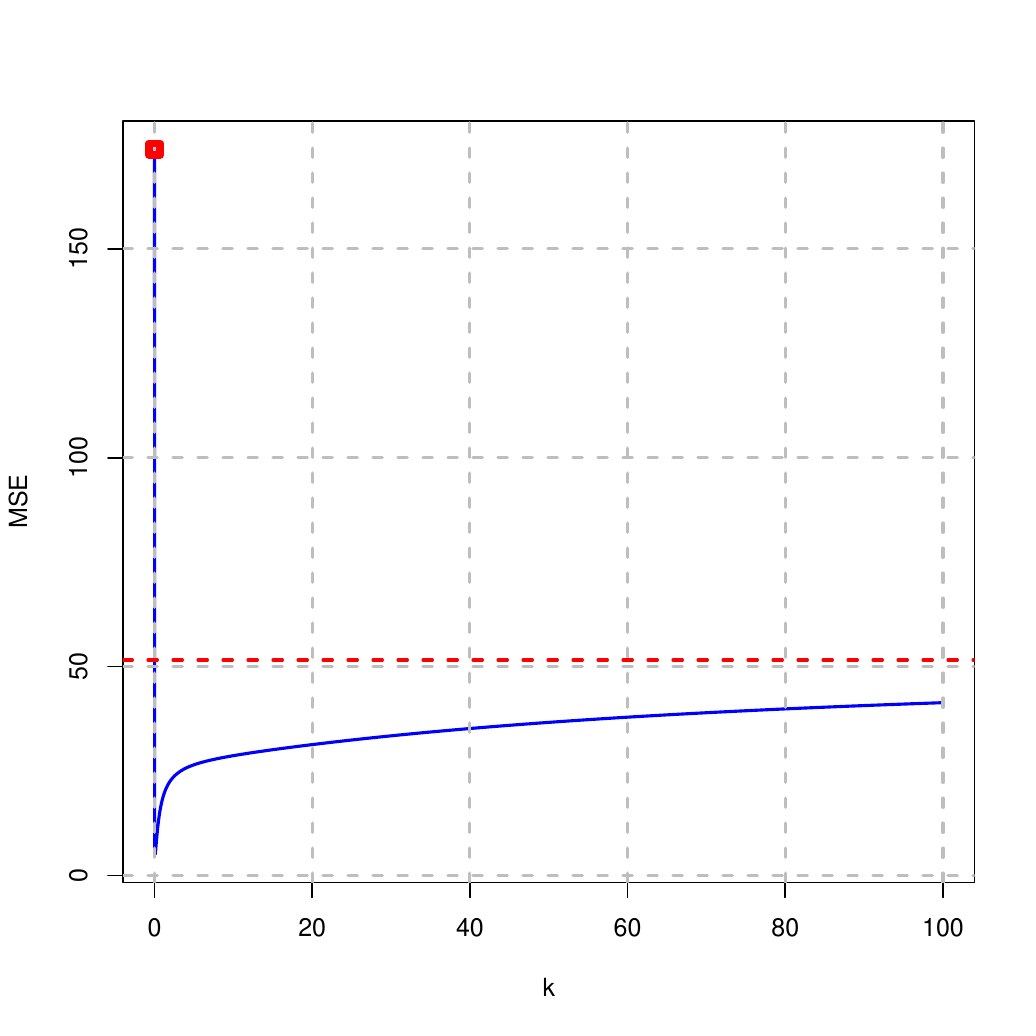}
  \includegraphics[width=8cm]{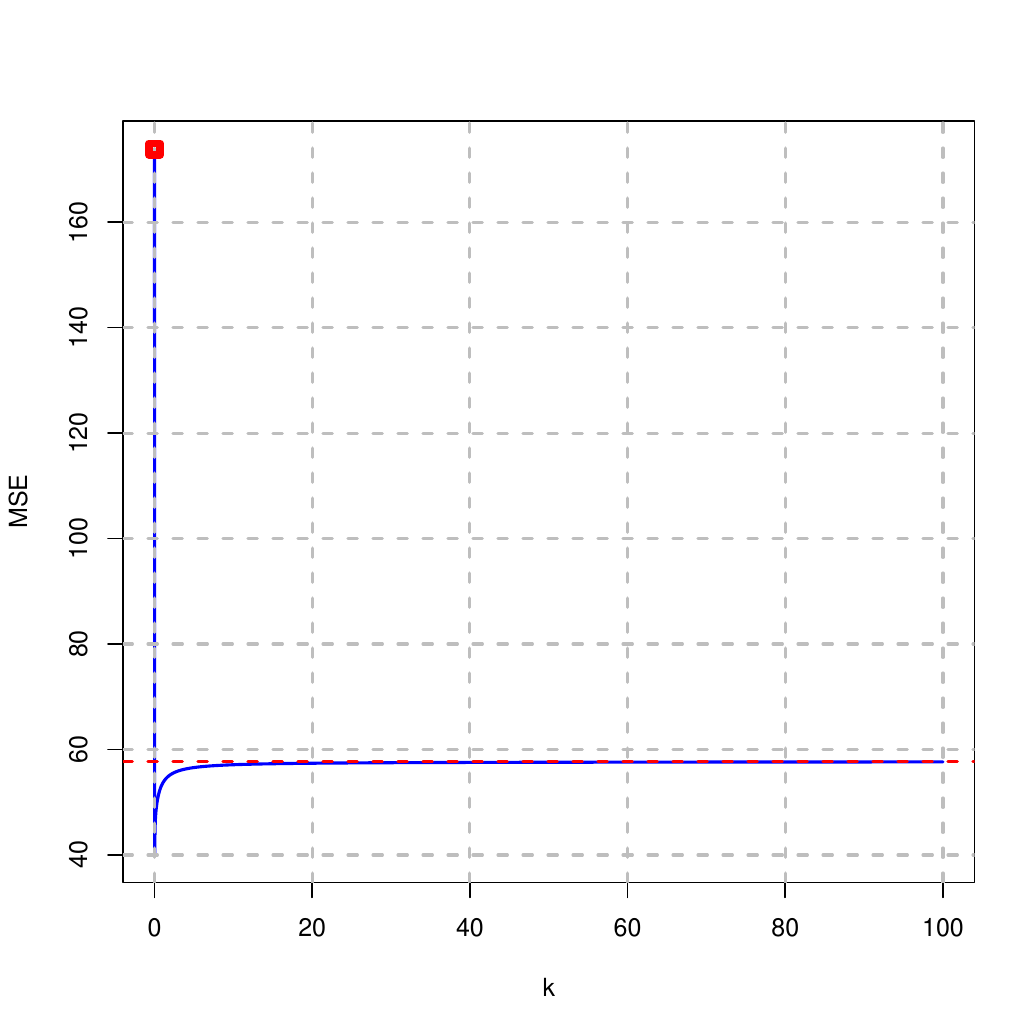}
  \caption{Trace of the mean squared errors of the penalized estimator ($h=1$, left) and ridge ($h=0$, right) of the model on bank credit in the United States for $k \in \{ 0, 0.01, 0.02, \dots, 100 \}$ (red dashed line represents the MSE asymptote)} \label{fig.example.4}
\end{figure}

\begin{figure}
  \centering
  \includegraphics[width=8cm]{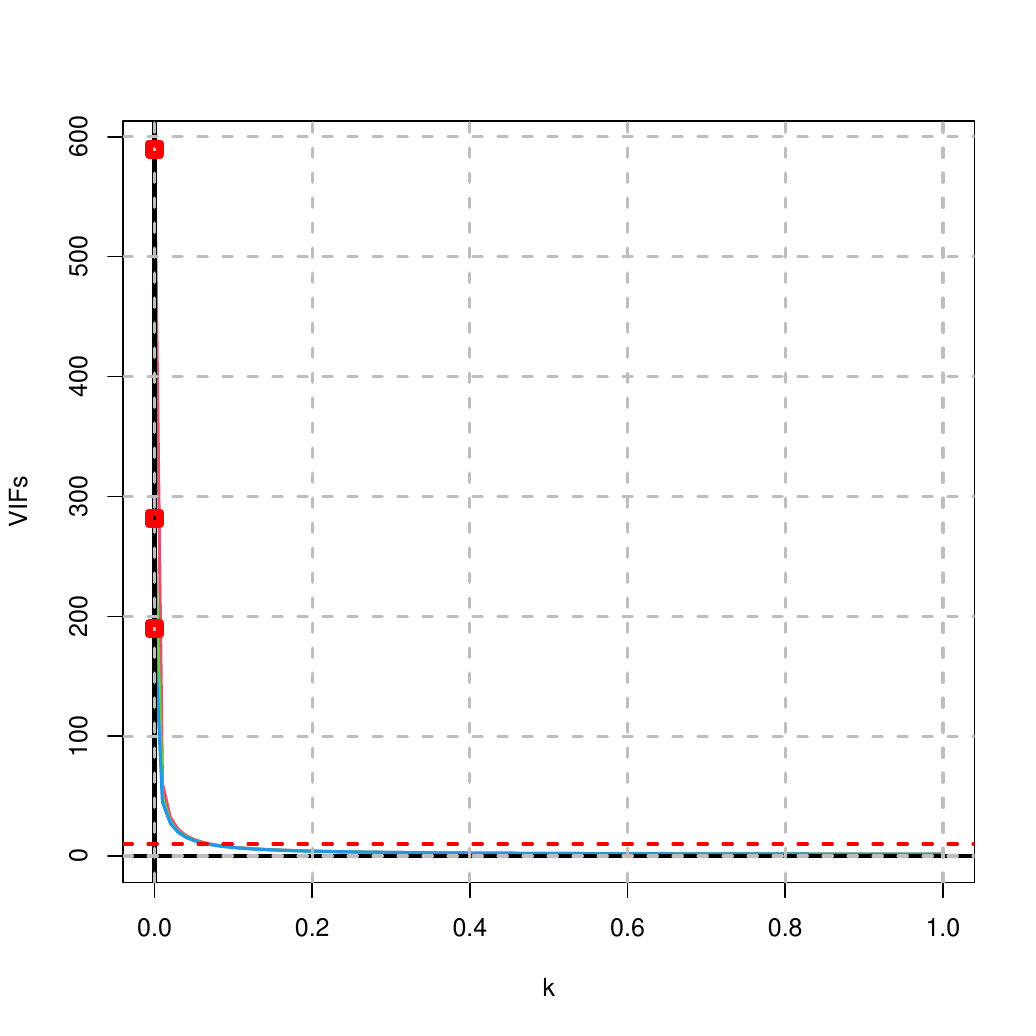}
  \includegraphics[width=8cm]{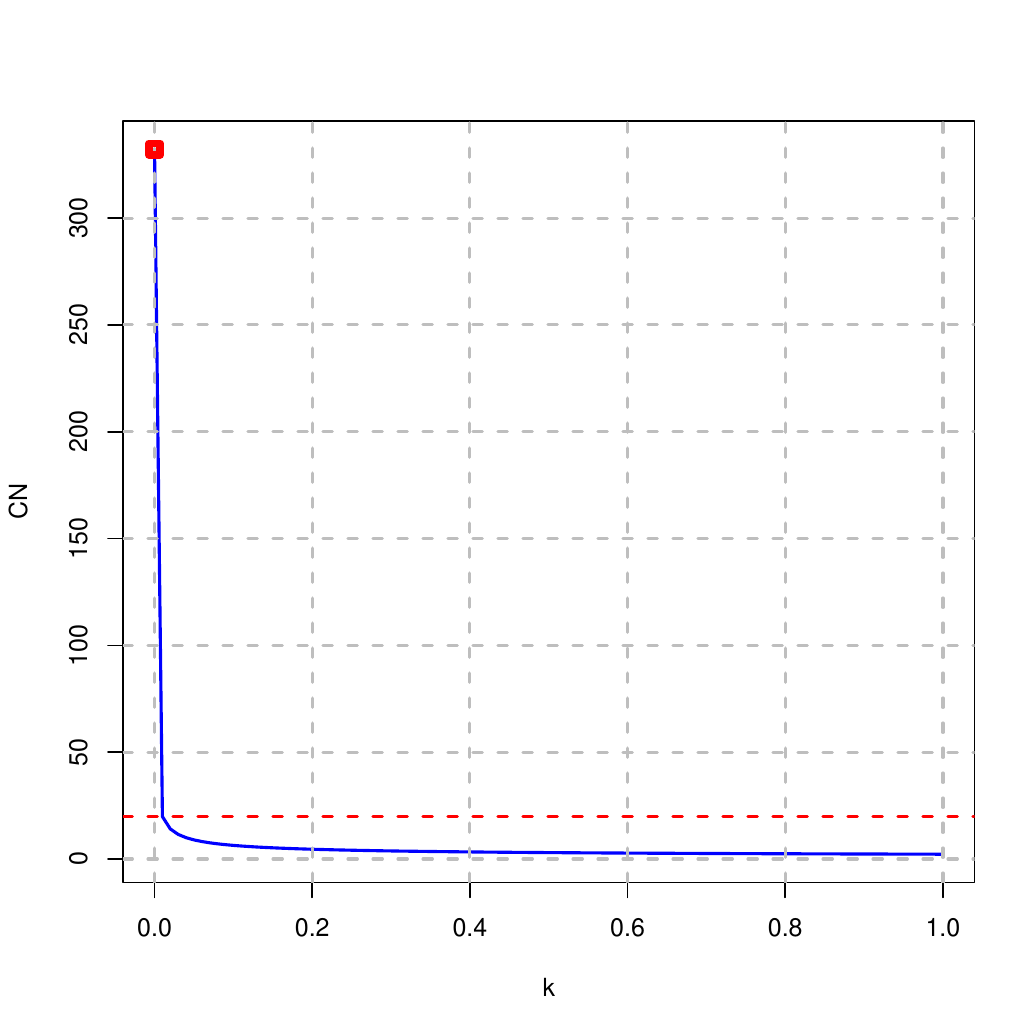}
  \caption{Trace of variance inflation factor (left) and condition number (right) of the model on bank credit in the United States for $k \in \{ 0, 0.01, 0.02, \dots, 1 \}$ (red dashed line represents the thresholds: 10 and 20, respectively.)} \label{fig.example.5}
\end{figure}

On the other hand, considering that$k \in \{ 0, 0.01, 0.02, \dots, 1 \}$, Figure \ref{fig.example.5} shows the traces of the variance inflation factor and the condition number. Note that in this case they coincide for the ridge and penalized estimator since the expanded matrix from which these measures are calculated coincide in both cases.

It is observed that in both cases they decrease towards their minimum values and that there are values of $k$ that are below the thresholds considered to be troubling.

\begin{figure}
  \centering
  \includegraphics[width=8cm]{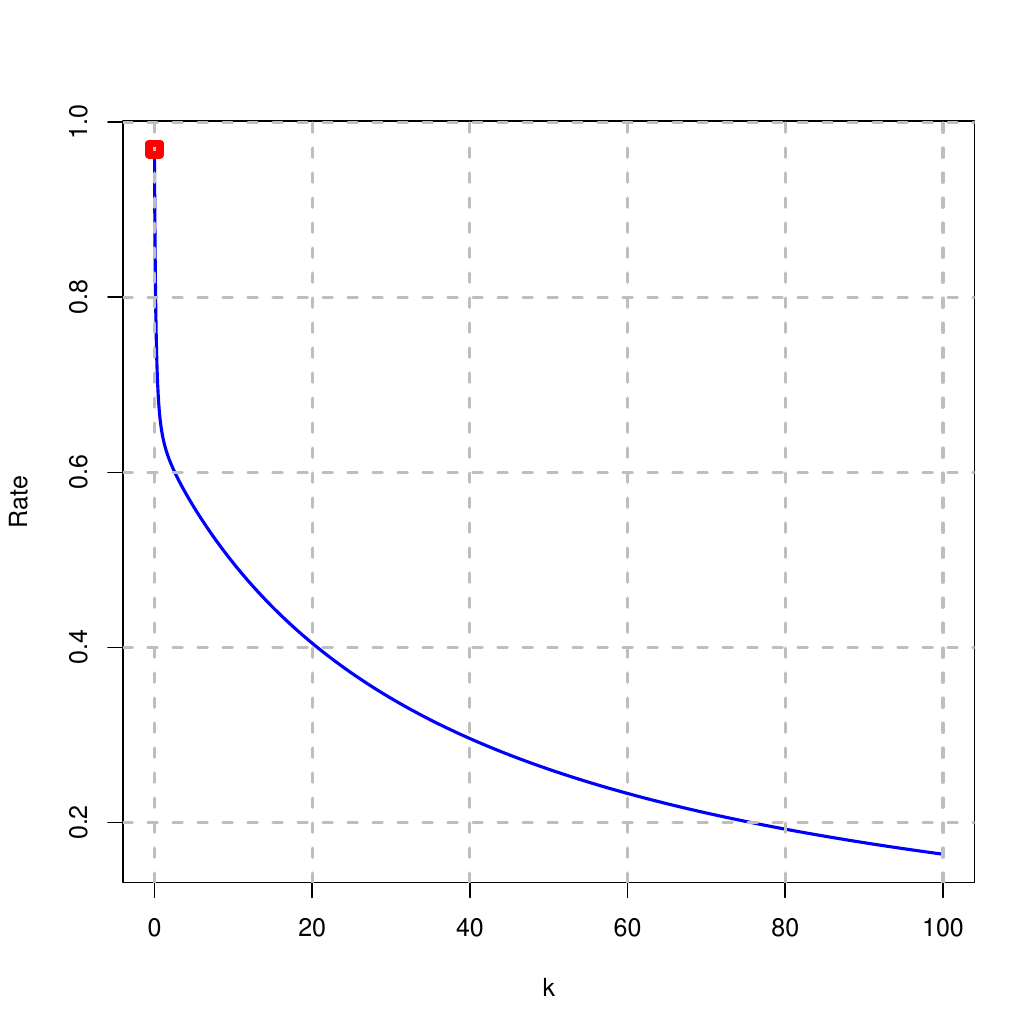}
  \caption{Trace of $|| \boldsymbol{\alpha} - \boldsymbol{\beta}(k,h) || / || \boldsymbol{\alpha} ||$ for the penalized estimator ($h=1$) of the model on bank credit in the United States for $k \in \{ 0, 0.01, 0.02, \dots, 100 \}$} \label{fig.example.6}
\end{figure}

Finally, considering once again that $k \in \{ 0, 0.01, 0.02, \dots, 100 \}$, Figure \ref{fig.example.6} shows the trace of $|| \boldsymbol{\alpha} - \boldsymbol{\beta}(k,h) || / || \boldsymbol{\alpha} ||$. It is observed that it is strictly monotonically decreasing, obtaining for the highest values of $k$ values less than 20\%.

\subsection{Choice of penalty parameter $k$ and analysis of the model}

As noted in the previous graphical representations, following the guidelines in subsection \ref{elegir_k}, it is possible to make a choice of $k$ that implies a variance inflation factor or condition number below the established thresholds or for a minimum mean square error. As a summary, Table \ref{tab:eleccion_k} shows the values of $k$ leading to the minimum square error and the first value of the variance inflation factor and condition number that falls below the thresholds traditionally set as troubling. Note that the penalized estimator has a smaller minimum square error than that of the ridge estimator.

\begin{table}
    \centering
    \begin{tabular}{ccc}
        \hline
        Criterion & Penalized estimator ($h=1$) & Ridge estimator ($h=0$) \\
        \hline
        minimum MSE & k=0.07 (5.131445) & k=0.02 (40.14453) \\
        NC lower than 30 & \multicolumn{2}{c}{k=0.01 (19.83053)} \\
        NC lower than 20 & \multicolumn{2}{c}{k=0.01 (19.83053)}  \\
        NC lower than 10 & \multicolumn{2}{c}{k=0.04 (9.966163)}  \\
        Maximum VIF lower than 10 & \multicolumn{2}{c}{k=0.08 (8.980033)}  \\
        \hline
    \end{tabular}
    \caption{Selection of the value of $k$ following the guidelines of the subsection \ref{elegir_k} (in parentheses the value of the measure for the value of $k$ under consideration)} \label{tab:eleccion_k}
\end{table}

For these values (including $k=0$), Table \ref{tab:est_cresta} shows the results obtained with the application of the ridge estimator while Table \ref{tab:est_penalizado} shows the results obtained with the application of the penalized estimator. Note that in this last table, the bootstrap approximations to the confidence interval obtained from the 0.025 and 0.975 percentiles are coded as type 1 and those obtained from the mean and quasi-variance as type 2.

\paragraph{Estimation of coefficients}

It is observed that the estimates of the ridge estimator decrease towards zero; contrary to those of the penalized estimator, although they are still far from the values fixed in the vector $\boldsymbol{\alpha}$.

On the other hand, inference using bootstrap methods from 10000 iterations indicates that the estimates obtained by the ridge estimator are not significantly different from zero for the $k$ values considered. However, this is not the case for the estimates obtained from the penalized estimator. Thus, for example, in the case where $k=0.08$:
\begin{itemize}
    \item The intercept is significantly different from zero with positive sign.
    \item The coefficient of personal consumption is significantly different from zero with a negative sign. That is, higher personal consumption is associated with a decrease in mortgage debt, which could be interpreted as higher consumption in situations where there is less mortgage debt. However, it should be noted that this estimated sign is the opposite  \footnote{
            Considering that $k=100$, it is obtained that $\widehat{\boldsymbol{\beta}}(h=1,k) = \left( 6.4705, 1.6741, 0.8247, -0.009 \right)^{t}$, all of them being coefficients significantly different from zero.In this case, the sign of the second coefficient has been corrected in accordance with the specifications of vector  $\boldsymbol{\alpha}$, although the fourth coefficient now has a negative sign.
            This situation is more understandable because the value set in the vector $\boldsymbol{\alpha}$ for this parameter, 0.00519, is susceptible to this possible change of sign.\\
            It is also obtained a root mean square error equal to 41.3982 (increases with respect to the minimum value although it is still lower than that of OLS), a condition number equal to 1.019 and maximum variance inflation factor equal to 1.000193. That is, the multicollinearity detection measures are very close to their minimum value, so the degree of multicollinearity is very low.  This issue is reflected in the fact that the average variance of the coefficient estimates when the independent variables are perturbed by 1\% is 0.1607\%.
        } of that fixed in the vector $\boldsymbol{\alpha}$.
    \item The coefficient of outstanding credit is significantly different from zero with a positive sign. That is, higher credit outstanding is associated with an increase in mortgage debt, which could be interpreted as meaning that there is higher credit outstanding when mortgage debt is higher.
\end{itemize}
Note that in this case we obtain a mean square error very close to the minimum value and much lower than that of OLS; at the same time we verify that the degree of existing approximate multicollinearity is considered not to be troubling since the number of condition and all the variance inflation factor are less than 10.

\paragraph{Goodness-of-fit}

Considering that for $k=0$ the results obtained coincide with those of OLS, it is observed that the value of the goodness-of-fit, 0.9878, does not coincide with the value of the coefficient of determination, 0.9235, provided in Table \ref{tab.OLS_Wissel}. This discrepancy is due, as discussed, to the fact that the proposed goodness-of-fit coincides with the coefficient of determination when $k=0$ only if the dependent variable has zero mean, a condition that in this case is not verified.

On the other hand, it is observed that the goodness-of-fit is very similar in the ridge and penalized estimator and, in addition, does not differ substantially from the value obtained in OLS.

\paragraph{Mean Square Error}

With respect to the mean square error, lower values are obtained in all cases in the case of the penalized estimator. This result is in line with the results shown in the Monte Carlo simulation of subsection \ref{simulatioMC}: the penalized estimator performs better in terms of MSE than the peak estimator when the troubling approximate multicollinearity is of the essential type.

This result is interesting as both the ridge and penalized estimators are biased estimators.

\paragraph{Existence of approximate multicollinearity}

As previously discussed, in both cases the same values are obtained for the variance inflation factor and the condition number. The condition number is less than 20 in all cases considered, while the maximum VIF is less than 10 only for $k=0.08$.

\paragraph{Numerical stability of the obtained estimates}

By perturbing the observations 10000 times by 1\% as indicated in the section \ref{stability} and calculating their effect on the coefficient estimates according to the expression (\ref{cuantification_perturb}), the average values shown in the tables are obtained.

It can be seen that initially small changes in the data result in an average variation in the coefficient estimates of 73.9726\%. This indicates a significant instability in the estimates obtained by OLS. This value decreases as the value of $k$ increases, the decrease being more pronounced in the case of the penalized estimator. Thus, for example, in the case where $k=0.08$ we have that a perturbation of 1\% in the independent variables implies an average variation of 7.6294\% in the coefficient estimates.

\begin{sidewaystable}
    \centering
    \begin{tabular}{ccccccccccc}
        \hline
            & \multicolumn{2}{c}{k=0} & \multicolumn{2}{c}{k=0.01} & \multicolumn{2}{c}{k=0.02} & \multicolumn{2}{c}{k=0.04} & \multicolumn{2}{c}{k=0.08} \\
            & Ridge & Bootstrap &  Ridge & Bootstrap &  Ridge & Bootstrap &  Ridge & Bootstrap &  Ridge & Bootstrap \\
        \hline
        $\beta_{1}$  & 5.4693 & 2.2482 & 1.0412 & 0.2857 & 0.2629 & -0.0943 & -0.21 & -0.3216 & -0.428 & -0.412 \\
        $\beta_{2}$  & -4.2524 & -2.466 & -2.4673 & -1.7351 & -2.1056 & -1.5493 & -1.814 & -1.379 & -1.5481 & -1.2016 \\
        $\beta_{3}$  & 3.1204 & 2.0841 & 2.5273 & 1.887 & 2.3586 & 1.789 & 2.1596 & 1.6495 & 1.8911 & 1.4441 \\
        $\beta_{4}$  & 0.0029 & 0.0028 & 0.0014 & 0.002 & 0.0013 & 0.002 & 0.0014 & 0.0021 & 0.0017 & 0.0025  \\
        DS for $\beta_{1}$  & 2947.0815 & 13.9162 & 902.4274 & 3.2142 & 532.7874 & 1.8279 & 292.8689 & 1.004 & 154.0945 & 0.5946 \\
        DS for $\beta_{2}$  & 1146.1535 & 6.1901 & 350.9741 & 2.167 & 207.2237 & 1.6559 & 113.9276 & 1.3197 & 59.9713 & 1.0715 \\
        DS for $\beta_{3}$  & 338.336 & 3.0476 & 103.6891 & 1.9856 & 61.3139 & 1.8011 & 33.8661 & 1.6106 & 18.0615 & 1.382 \\
        DS for $\beta_{4}$  & 1.0981 & 0.0055 & 0.3364 & 0.0032 & 0.1987 & 0.0031 & 0.1095 & 0.0029 & 0.0579 & 0.0025 \\
        GoF  & 0.9878 & 0.991 & 0.9877 & 0.9905 & 0.9876 & 0.9903 & 0.9876 & 0.9902 & 0.9874 & 0.99 \\
        \hline
        BIT1 for $\beta_{1}$  & \multicolumn{2}{c}{(-23.3118, 28.5705)} & \multicolumn{2}{c}{(-5.4661, 6.7827)} & \multicolumn{2}{c}{(-3.4506, 3.648)} & \multicolumn{2}{c}{(-2.2583, 1.7004)} & \multicolumn{2}{c}{(-1.6039, 0.6764)} \\
        BIT1 for $\beta_{2}$  & \multicolumn{2}{c}{(-13.1453, 9.7675)} & \multicolumn{2}{c}{(-5.1669, 2.8114)} & \multicolumn{2}{c}{(-4.0082, 1.9754)} & \multicolumn{2}{c}{(-3.2576, 1.4472)} & \multicolumn{2}{c}{(-2.7288, 1.0429)} \\
        BIT1 for $\beta_{3}$  & \multicolumn{2}{c}{(-4.5941, 6.4143)} & \multicolumn{2}{c}{(-2.4435, 4.7682)} & \multicolumn{2}{c}{(-2.1057, 4.4248)} & \multicolumn{2}{c}{(-1.7724, 4.0229)} & \multicolumn{2}{c}{(-1.4312, 3.5192)} \\
        BIT1 for $\beta_{4}$  & \multicolumn{2}{c}{(-0.0082, 0.0135)} & \multicolumn{2}{c}{(-0.004, 0.0081)} & \multicolumn{2}{c}{(-0.0035, 0.0081)} & \multicolumn{2}{c}{(-0.0028, 0.0079)} & \multicolumn{2}{c}{(-0.0018, 0.0074)} \\
        BIT1 for GoF   & \multicolumn{2}{c}{(0.9852, 0.998)} & \multicolumn{2}{c}{(0.9848, 0.9974)} & \multicolumn{2}{c}{(0.9845, 0.9972)} & \multicolumn{2}{c}{(0.9843, 0.997)} & \multicolumn{2}{c}{(0.9841, 0.9968)} \\
        BIT2 for $\beta_{1}$  & \multicolumn{2}{c}{(-25.0275, 29.5239)} & \multicolumn{2}{c}{(-6.0141, 6.5856)} & \multicolumn{2}{c}{(-3.6769, 3.4884)} & \multicolumn{2}{c}{(-2.2895, 1.6464)} & \multicolumn{2}{c}{(-1.5774, 0.7535)} \\
        BIT2 for $\beta_{2}$  & \multicolumn{2}{c}{(-14.5985, 9.6666)} & \multicolumn{2}{c}{(-5.9824, 2.5121)} & \multicolumn{2}{c}{(-4.7948, 1.6962)} & \multicolumn{2}{c}{(-3.9655, 1.2076)} & \multicolumn{2}{c}{(-3.3018, 0.8986)} \\
        BIT2 for $\beta_{3}$  & \multicolumn{2}{c}{(-3.8892, 8.0574)} & \multicolumn{2}{c}{(-2.0049, 5.7788)} & \multicolumn{2}{c}{(-1.7411, 5.3191)} & \multicolumn{2}{c}{(-1.5072, 4.8063)} & \multicolumn{2}{c}{(-1.2646, 4.1529)} \\
        BIT2 for $\beta_{4}$  & \multicolumn{2}{c}{(-0.0081, 0.0136)} & \multicolumn{2}{c}{(-0.0043, 0.0083)} & \multicolumn{2}{c}{(-0.004, 0.008)} & \multicolumn{2}{c}{(-0.0035, 0.0077)} & \multicolumn{2}{c}{(-0.0024, 0.0074)} \\
        BIT2 for GoF   & \multicolumn{2}{c}{(0.9847, 0.9973)} & \multicolumn{2}{c}{(0.9843, 0.9967)} & \multicolumn{2}{c}{(0.9841, 0.9966)} & \multicolumn{2}{c}{(0.9839, 0.9965)} & \multicolumn{2}{c}{(0.9837, 0.9963)} \\
        \hline
        MSE  & \multicolumn{2}{c}{199.9497} & \multicolumn{2}{c}{44.3756} & \multicolumn{2}{c}{41.3225} & \multicolumn{2}{c}{43.3614} & \multicolumn{2}{c}{45.9795} \\
        CN  & \multicolumn{2}{c}{332.3} & \multicolumn{2}{c}{19.8305} & \multicolumn{2}{c}{14.0525} & \multicolumn{2}{c}{9.9662} & \multicolumn{2}{c}{7.0841} \\
        VIF(1,k)  & \multicolumn{2}{c}{589.754} & \multicolumn{2}{c}{60.3983} & \multicolumn{2}{c}{32.2129} & \multicolumn{2}{c}{16.9411} & \multicolumn{2}{c}{8.98} \\
        VIF(1,k)  & \multicolumn{2}{c}{281.8862} & \multicolumn{2}{c}{48.7101} & \multicolumn{2}{c}{28.4015} & \multicolumn{2}{c}{15.8204} & \multicolumn{2}{c}{8.6686} \\
        VIF(1,k)  & \multicolumn{2}{c}{189.4874} & \multicolumn{2}{c}{45.2022} & \multicolumn{2}{c}{27.2576} & \multicolumn{2}{c}{15.4841} & \multicolumn{2}{c}{8.5752} \\
        \hline
        Stability  & \multicolumn{2}{c}{73.9726} & \multicolumn{2}{c}{50.4276} & \multicolumn{2}{c}{36.4043} & \multicolumn{2}{c}{24.5117} & \multicolumn{2}{c}{17.0305} \\
        \hline
    \end{tabular}
    \caption{Analysis of the model on bank credit in the United States using the ridge estimator for $k=0.01, 0.02, 0.04, 0.08$ (BIT1 = Bootstrap interval type 1; BIT2 = Bootstrap interval type 2; DS = Desviation Standard)} \label{tab:est_cresta}
\end{sidewaystable}

\begin{sidewaystable}
    \centering
    \begin{tabular}{ccccccccccc}
        \hline
            & \multicolumn{2}{c}{k=0} & \multicolumn{2}{c}{k=0.01} & \multicolumn{2}{c}{k=0.04} & \multicolumn{2}{c}{k=0.07} & \multicolumn{2}{c}{k=0.08} \\
            & Penalized & Bootstrap &  Penalized & Bootstrap &  Penalized & Bootstrap &  Penalized & Bootstrap &  Penalized & Bootstrap \\
        \hline
        $\beta_{1}$  & 5.4693 & 2.2482 & 4.585 & 4.0199 & 4.4378 & 4.3899 & 4.4938 & 4.5378 & 4.5171 & 4.5763  \\
        $\beta_{2}$  & -4.2524 & -2.466 & -3.8019 & -3.1602 & -3.4622 & -3.0343 & -3.247 & -2.8617 & -3.1853 & -2.8093 \\
        $\beta_{3}$  & 3.1204 & 2.0841 & 2.8777 & 2.2792 & 2.4876 & 1.9617 & 2.182 & 1.6749 & 2.0917 & 1.5909 \\
        $\beta_{4}$  & 0.0029 & 0.0028 & 0.0028 & 0.0035 & 0.0035 & 0.0043 & 0.0041 & 0.0049 & 0.0042 & 0.0051 \\
        DS for $\beta_{1}$  & 2947.0815 & 13.9162 & 902.4274 & 3.1633 & 292.8689 & 0.9857 & 174.8012 & 0.6339 & 154.0945 & 0.5787 \\
        DS for $\beta_{2}$  & 1146.1535 & 6.1901 & 350.9741 & 2.1118 & 113.9276 & 1.3152 & 68.0214 & 1.1398 & 59.9713 & 1.1021 \\
        DS for $\beta_{3}$  & 338.336 & 3.0476 & 103.6891 & 1.9337 & 33.8661 & 1.6009 & 20.4131 & 1.4532 & 18.0615 & 1.4139 \\
        DS for $\beta_{4}$  & 1.0981 & 0.0055 & 0.3364 & 0.0032 & 0.1095 & 0.0028 & 0.0656 & 0.0026 & 0.0579 & 0.0025 \\
        GoF  & 0.9878 & 0.991 & 0.9878 & 0.9905 & 0.9877 & 0.9901 & 0.9875 & 0.9899 & 0.9874 & 0.9898 \\
        \hline
        BIT1 for $\beta_{1}$  & \multicolumn{2}{c}{(-23.3118, 28.5705)} & \multicolumn{2}{c}{(-1.615, 10.4396)} & \multicolumn{2}{c}{(2.5116, 6.3897)} & \multicolumn{2}{c}{(3.2787, 5.7392)} & \multicolumn{2}{c}{(3.4191, 5.6548)}  \\
        BIT1 for $\beta_{2}$  & \multicolumn{2}{c}{(-13.1453, 9.7675)} & \multicolumn{2}{c}{(-6.4747, 1.3367)} & \multicolumn{2}{c}{(-4.8868, -0.193)} & \multicolumn{2}{c}{(-4.4573, -0.4378)} & \multicolumn{2}{c}{(-4.3505, -0.4793)} \\
        BIT1 for $\beta_{3}$  & \multicolumn{2}{c}{(-4.5941, 6.4143)} & \multicolumn{2}{c}{(-1.9709, 5.1374)} & \multicolumn{2}{c}{(-1.4629, 4.3341)} & \multicolumn{2}{c}{(-1.3739, 3.8329)} & \multicolumn{2}{c}{(-1.3547, 3.6979)} \\
        BIT1 for $\beta_{4}$  & \multicolumn{2}{c}{(-0.0082, 0.0135)} & \multicolumn{2}{c}{(-0.0025, 0.0096)} & \multicolumn{2}{c}{(-0.0006, 0.01)} & \multicolumn{2}{c}{(0.00005, 0.0101)} & \multicolumn{2}{c}{(0.0008, 0.0101)} \\
        BIT1 for GoF   & \multicolumn{2}{c}{(0.9852, 0.998)} & \multicolumn{2}{c}{(0.9849, 0.9969)} & \multicolumn{2}{c}{(0.9845, 0.9964)} & \multicolumn{2}{c}{(0.9842, 0.9963)} & \multicolumn{2}{c}{(0.9841, 0.9962)} \\
        BIT2 for $\beta_{1}$  & \multicolumn{2}{c}{(-25.0275, 29.5239)} & \multicolumn{2}{c}{(-2.1802, 10.2199)} & \multicolumn{2}{c}{(2.4578, 6.3219)} & \multicolumn{2}{c}{(3.2954, 5.7803)} & \multicolumn{2}{c}{(3.442, 5.7106)} \\
        BIT2 for $\beta_{2}$  & \multicolumn{2}{c}{(-14.5985, 9.6666)} & \multicolumn{2}{c}{(-7.2993, 0.9789)} & \multicolumn{2}{c}{(-5.612, -0.4566)} & \multicolumn{2}{c}{(-5.0956, -0.6277)} & \multicolumn{2}{c}{(-4.9695, -0.6491)} \\
        BIT2 for $\beta_{3}$  & \multicolumn{2}{c}{(-3.8892, 8.0574)} & \multicolumn{2}{c}{(-1.5109, 6.0694)} & \multicolumn{2}{c}{(-1.1761, 5.0996)} & \multicolumn{2}{c}{(-1.1733, 4.5232)} & \multicolumn{2}{c}{(-1.1804, 4.3622)} \\
        BIT2 for $\beta_{4}$  & \multicolumn{2}{c}{(-0.0081, 0.0136)} & \multicolumn{2}{c}{(-0.0027, 0.0097)} & \multicolumn{2}{c}{(-0.0013, 0.0099)} & \multicolumn{2}{c}{(-0.0002, 0.01)} & \multicolumn{2}{c}{(0.0001, 0.0101)} \\
        BIT2 for GoF   & \multicolumn{2}{c}{(0.9847, 0.9973)} & \multicolumn{2}{c}{(0.9845, 0.9965)} & \multicolumn{2}{c}{(0.9842, 0.9961)} & \multicolumn{2}{c}{(0.9839, 0.9959)} & \multicolumn{2}{c}{(0.9838, 0.9958)} \\
        \hline
        MSE  & \multicolumn{2}{c}{199.9497} & \multicolumn{2}{c}{22.2729} & \multicolumn{2}{c}{6.3279} & \multicolumn{2}{c}{5.4749} & \multicolumn{2}{c}{5.4808} \\
        CN  & \multicolumn{2}{c}{332.3} & \multicolumn{2}{c}{19.8305} & \multicolumn{2}{c}{9.9662} & \multicolumn{2}{c}{7.5635} & \multicolumn{2}{c}{7.0841} \\
        VIF(1,k)  & \multicolumn{2}{c}{589.754} & \multicolumn{2}{c}{60.3983} & \multicolumn{2}{c}{16.9411} & \multicolumn{2}{c}{10.131} & \multicolumn{2}{c}{8.98} \\
        VIF(2,k)  & \multicolumn{2}{c}{281.8862} & \multicolumn{2}{c}{48.7101} & \multicolumn{2}{c}{15.8204} & \multicolumn{2}{c}{9.7315} & \multicolumn{2}{c}{8.6686} \\
        VIF(3,k)  & \multicolumn{2}{c}{189.4874} & \multicolumn{2}{c}{45.2022} & \multicolumn{2}{c}{15.4841} & \multicolumn{2}{c}{9.6116} & \multicolumn{2}{c}{8.5752} \\
        \hline
        Stability  & \multicolumn{2}{c}{73.9726} & \multicolumn{2}{c}{29.8586} & \multicolumn{2}{c}{12.0844} & \multicolumn{2}{c}{8.3133} & \multicolumn{2}{c}{7.6294} \\
        \hline
    \end{tabular}
    \caption{Analysis of the model on bank credit in the United States using the penalized estimator for $k=0.01, 0.04, 0.07, 0.08$ (BIT1 = Bootstrap interval type 1; BIT2 = Bootstrap interval type 2; DS = Desviation Standard)} \label{tab:est_penalizado}
\end{sidewaystable}

\section{Conclusions} \label{conc.SCRP}

This paper proposes the estimation of the coefficients of a multiple linear regression model by penalizing the function to be minimized with the objective of getting the estimates to be as close as possible to the relationships given by the corresponding simple linear regressions, since then the relationships of each of the independent variables to the dependent variable are obtained if the rest of the
explanatory variables vary as would be expected.

The main characteristics of the proposed penalized estimator (trace and norm of the estimator, variance-covariance matrix, goodness-of-fit and mean square error) and its implications in the detection and treatment of multicollinearity have been thoroughly analyzed. This analysis has served to propose different possibilities for choosing the penalty parameter, so that for a fixed value of this parameter, inference can be performed by  bootstrap methodology. The fact that the objective is that the estimates converge to the estimates of the simple regressions means that the bootstrap inference proposed can provide estimates of the coefficients of the independent variables significantly different from zero, contrary to what occurs in the ridge estimator, since in this case it is established that the estimates converge to zero.

It has also been found that the ridge estimator of Hoerl and Kennard \cite{HKa, HKb} is a particular case of the proposed estimator. Due to this relationship, several characteristics coincide in both estimators, such as the variance-covariance matrix, the calculation of the variance inflation factor (VIF) and the condition number (CN). This makes it possible that studies on these measures in the ridge estimator will be directly applicable to the penalized estimator obtained.

It is also important to note that this estimator mitigates the numerical instability that can occur when estimating by OLS a multiple linear regression model where the approximate multicollinearity is troubling, so its application is interesting when this consequence of multicollinearity is present in the econometric model analyzed.

On the other hand, in the works of García et al. \cite{Gracia2015}, Salmerón  et al. \cite{Salmeron2017} and Rodríguez et al. \cite{Rodriguez2022,Rodriguezetal2021} where the VIF, CN, goodness of fit and Stewart index, are respectively extended (see \cite{Stewart1987}) to ridge estimator, it is stated that \textit{the transformation of the data is not optional but is required for the correct application of these measures in ridge regression}. This issue has already been addressed by Marquardt in his work entitled ``You should standardize the predictor variables in your regression models'' (\cite{Marquardt1980}).
While Obenchain \cite{Obenchain1975} indicated that when working with the ridge estimator \textit{the anaylisis should be completely redone if any transformation of regressors is adopted}.
In this sense, this work establishes that for the penalized estimator it is not necessary to perform any transformation of the data except to consider the standardization of the data to calculate the VIF and the unit length to calculate the CN.

Finally, future lines of work include the combination with the LASSO estimator to emulate the elastic net methodology and perform variable selection,  analyze the possibilities of the proposed estimator when there are outliers in the data and compare its performance in terms of MSE with the proposals from Wang et al. \cite{Wangetal2021} and Lukman et al. \cite{Lukman2024}. A future goal is also the creation of a package in R \cite{RCoreTeam} from the code available in Github that allows the application of the methodology presented to any user in the teaching, research and/or business environment.

\appendix

\section{Results of interest}
    \label{apen.operaciones}

\subsection{Norm of the penalized estimator}
    \label{norma.penalizada}

    Taking into account the decomposition $\mathbf{X}^{t} \mathbf{X} = \boldsymbol{\Gamma} \mathbf{D}_{\lambda_{i}} \boldsymbol{\Gamma}^{t}$ donde $\boldsymbol{\Gamma}_{p \times p}$ is an orthogonal matrix ($\boldsymbol{\Gamma}^{t} \boldsymbol{\Gamma} = \mathbf{I} = \boldsymbol{\Gamma} \boldsymbol{\Gamma}^{t}$) which contains the eigenvectors of $\mathbf{X}^{t} \mathbf{X}$ and $\mathbf{D}_{\lambda_{i}}$ is a diagonal matrix ($\mathbf{D}_{\lambda_{i}} = diag \left( \lambda_{1}, \dots, \lambda_{p} \right)$) of dimension $p \times p$ containing its eigenvalues (which are real and positive), it is verified that:
    \begin{equation}
        \label{Z.apen}
        \mathbf{Z}(k) = \left( \mathbf{X}^{t} \mathbf{X} + k \mathbf{I} \right)^{-1} = \boldsymbol{\Gamma} \mathbf{D}_{\frac{1}{\lambda_{i}+k}} \boldsymbol{\Gamma}^{t}.
    \end{equation}

    Denoting $\boldsymbol{\Psi}_{p \times 1} = \mathbf{X}^{t} \mathbf{y}$ in expression (\ref{est.SCRP.1}), it is obtained that:
    \begin{eqnarray*}
        || \widehat{\boldsymbol{\beta}}(k, h) || &=& \widehat{\boldsymbol{\beta}}(k, h)^{t} \widehat{\boldsymbol{\beta}}(k, h)
            = \left( \boldsymbol{\Psi} + k \cdot h \cdot \boldsymbol{\alpha} \right)^{t} \cdot \mathbf{Z}(k)^{t} \mathbf{Z}(k) \cdot \left( \boldsymbol{\Psi} + k \cdot h \cdot \boldsymbol{\alpha} \right)  \\
            &=& \boldsymbol{\Psi}^{t} \boldsymbol{\Gamma} \mathbf{D}_{\frac{1}{(\lambda_{i}+k)^{2}}} \boldsymbol{\Gamma}^{t} \boldsymbol{\Psi}
                + 2 \cdot k \cdot h \cdot \boldsymbol{\Psi}^{t} \boldsymbol{\Gamma} \mathbf{D}_{\frac{1}{(\lambda_{i}+k)^{2}}} \boldsymbol{\Gamma}^{t} \boldsymbol{\alpha}
                + k^{2} \cdot h^{2}  \cdot \boldsymbol{\alpha}^{t} \boldsymbol{\Gamma} \mathbf{D}_{\frac{1}{(\lambda_{i}+k)^{2}}} \boldsymbol{\Gamma}^{t} \boldsymbol{\alpha} \\
            &=& \mathbf{a}^{t} \mathbf{D}_{\frac{1}{(\lambda_{i}+k)^{2}}} \mathbf{a}
                + 2 \cdot h \cdot \mathbf{a}^{t} \mathbf{D}_{\frac{k}{(\lambda_{i}+k)^{2}}} \mathbf{b}
                + h^{2} \cdot \mathbf{b}^{t} \mathbf{D}_{\frac{k^{2}}{(\lambda_{i}+k)^{2}}} \mathbf{b} \\
            &=& \sum \limits_{i=1}^{p} \frac{a_{i}^{2}}{(\lambda_{i}+k)^{2}} + 2 \cdot h \cdot \sum \limits_{i=1}^{p} \frac{k a_{i} b_{i}}{(\lambda_{i}+k)^{2}} + h^{2} \cdot \sum \limits_{i=1}^{p} \frac{k^{2} b_{i}^{2}}{(\lambda_{i}+k)^{2}},
    \end{eqnarray*}
    where it was considered that $\mathbf{a}_{p \times 1} = \boldsymbol{\Gamma}^{t} \boldsymbol{\Psi}$ and $\mathbf{b}_{p \times 1} = \boldsymbol{\Gamma}^{t} \boldsymbol{\alpha}$.

    Deriving with respect to $k$:
    \begin{eqnarray}
        \frac{\partial}{\partial k} \left( \sum \limits_{i=1}^{p} \frac{a_{i}^{2}}{(\lambda_{i}+k)^{2}} \right) &=& - \sum \limits_{i=1}^{p} \frac{2 a_{i}^{2}}{(\lambda_{i}+k)^{3}}, \label{deriv1} \\
        \frac{\partial}{\partial k} \left( \sum \limits_{i=1}^{p} \frac{k a_{i} b_{i}}{(\lambda_{i}+k)^{2}} \right) &=& \sum \limits_{i=1}^{p} \frac{(\lambda_{i}-k) a_{i} b_{i}}{(\lambda_{i}+k)^{3}}, \label{deriv2} \\
        \frac{\partial}{\partial k} \left( \sum \limits_{i=1}^{p} \frac{k^{2} b_{i}^{2}}{(\lambda_{i}+k)^{2}} \right) &=&  \sum \limits_{i=1}^{p} \frac{2 k \lambda_{i} b_{i}^{2}}{(\lambda_{i}+k)^{3}}, \label{deriv3}
    \end{eqnarray}
    the first term is decreasing when $k$ increases while the third term is increasing since $k, \lambda_{i}>0$ and then its derivatives are, expressions  (\ref{deriv1}) and (\ref{deriv3}) respectively, negatives and positives. At the same time the sign of the derivative of the second term (expression (\ref{deriv2})) is undefined since, on the one hand, it depends on the sign of $\lambda_{i}-k$ and, on the other hand, on the sign of $a_{i}b_{i}$.

    In short, when $h \not=0$ cannot be concluded about the monotony of the $|| \widehat{\boldsymbol{\beta}}(k, h) ||$, although it is certain that there exists a value of $k$ that stabilizes this norm, since when $k \rightarrow +\infty$:
    $$|| \widehat{\boldsymbol{\beta}}(k, h) || \rightarrow h^{2} \cdot \sum \limits_{i=1}^{p} b_{i}^{2} = h^{2} \cdot \mathbf{b}^{t} \mathbf{b} = h^{2} \cdot \boldsymbol{\alpha}^{t} \boldsymbol{\Gamma}^{t} \boldsymbol{\Gamma} \boldsymbol{\alpha} = h^{2} \cdot \boldsymbol{\alpha}^{t} \boldsymbol{\alpha} =  h^{2} \cdot || \boldsymbol{\alpha}||.$$

    Finally, when $h=0$ (the penalized estimator coincides with the ridge estimator), it is clear that $|| \widehat{\boldsymbol{\beta}}(k, h) ||$ is decreasing in $k$.

\subsection{Optimality of the variance-covariance matrix}
    \label{varianza.optima}

    Starting with the expression (\ref{Z.apen}), it follows that $\mathbf{Z}(k) \mathbf{X}^{t} \mathbf{X} \mathbf{Z}(k) = \boldsymbol{\Gamma} \mathbf{D}_{\frac{\lambda{i}}{(\lambda_{i} + k)^{2}}} \boldsymbol{\Gamma}^{t}$. Since $\left( \mathbf{X}^{t} \mathbf{X} \right)^{-1} = \boldsymbol{\Gamma} \mathbf{D}_{\frac{1}{\lambda_{i}}} \boldsymbol{\Gamma}^{t}$, it is clear that:
    $$var \left( \widehat{\boldsymbol{\beta}}(k, h) \right) - var \left( \widehat{\boldsymbol{\beta}} \right) = \sigma^{2} \cdot \boldsymbol{\Gamma} \mathbf{D}_{\frac{\lambda_{i}^{2} - (\lambda_{i}+k)^{2}}{\lambda_{i} (\lambda_{i}+k)^{2}}} \boldsymbol{\Gamma}^{t}.$$

    Considering $\mathbf{c}_{p \times 1}$ a non-zero vector and taking into account that $k>0$ and $\lambda_{i} > 0$ for $i=1,\dots,p$, it is verified that:
    $$\mathbf{c} \boldsymbol{\Gamma} \mathbf{D}_{\frac{\lambda_{i}^{2} - (\lambda_{i}+k)^{2}}{\lambda_{i} (\lambda_{i}+k)^{2}}} \boldsymbol{\Gamma}^{t}  \mathbf{c} = \mathbf{d} \mathbf{D}_{\frac{\lambda_{i}^{2} - (\lambda_{i}+k)^{2}}{\lambda_{i} (\lambda_{i}+k)^{2}}} \mathbf{d} = \sum \limits_{i=1}^{p} \frac{\lambda_{i}^{2} - (\lambda_{i}+k)^{2}}{\lambda_{i} (\lambda_{i}+k)^{2}} d_{i}^{2} = - \sum \limits_{i=1}^{p} \frac{k^{2} + 2 \lambda_{i} k}{\lambda_{i} (\lambda_{i}+k)^{2}} d_{i}^{2} < 0,$$
    where $\mathbf{d}_{p \times 1} = \boldsymbol{\Gamma}^{t} \mathbf{c}$. Therefore, $\boldsymbol{\Gamma} \mathbf{D}_{\frac{\lambda_{i}^{2} - (\lambda_{i}+k)^{2}}{\lambda_{i} (\lambda_{i}+k)^{2}}} \boldsymbol{\Gamma}^{t}$ es a negative definite matrix and, since $\sigma^{2} >0$, it is verified:
    $$var \left( \widehat{\boldsymbol{\beta}}(k, h) \right) - var \left( \widehat{\boldsymbol{\beta}} \right) < 0 \rightarrow var \left( \widehat{\boldsymbol{\beta}}(k, h) \right) < var \left( \widehat{\boldsymbol{\beta}} \right).$$

\subsection{Estimation errors}
    \label{ba.errores}

    Following the development of Rodríguez et al. \cite{Rodriguez2022} (Proposition 1), taking into account that the expression (\ref{est.SCRP.1}) comes from the system of normal equations:
    $$\left( \mathbf{X}^{t} \mathbf{X} + k \cdot \mathbf{I} \right) \cdot \widehat{\boldsymbol{\beta}}(k, h) = \mathbf{X}^{t} \mathbf{y} + k \cdot h \cdot \boldsymbol{\alpha},$$
    it must be verified that the first row of $\mathbf{X}^{t} \mathbf{X} + k \cdot \mathbf{I}$ , multiplied by  $\widehat{\boldsymbol{\beta}}(k, h)$:
    \begin{eqnarray*}
        \left( n+k, \sum \limits_{j=1}^{p} X_{2j}, \cdots, \sum \limits_{j=1}^{p} X_{pj} \right) \cdot \left(
        \begin{array}{c}
            \widehat{\beta}_{1}(k, h) \\
            \widehat{\beta}_{2}(k, h) \\
            \vdots \\
            \widehat{\beta}_{p}(k, h) \\
        \end{array} \right) &=& (n+k) \cdot \widehat{\beta}_{1}(k, h) + \widehat{\beta}_{2}(k, h) \cdot \sum \limits_{j=1}^{p} X_{2j} \\
        & & + \cdots  + \widehat{\beta}_{p}(k, h) \cdot \sum \limits_{j=1}^{p} X_{pj},
    \end{eqnarray*}
    has to be equal to the first element of $\mathbf{X}^{t} \mathbf{y} + k \cdot h \cdot \boldsymbol{\alpha}$, that is to say to $\sum \limits_{j=1}^{n} y_{j} + k \cdot h \cdot \overline{\mathbf{y}}$. Consequently:
    $$\sum \limits_{j=1}^{n} y_{j} + k \cdot h \cdot \overline{\mathbf{y}} = (n+k) \cdot \widehat{\beta}_{1}(k, h) + \widehat{\beta}_{2}(k, h) \cdot \sum \limits_{j=1}^{p} X_{2j} + \cdots  + \widehat{\beta}_{p}(k, h) \cdot \sum \limits_{j=1}^{p} X_{pj},$$
    or equivalently:
    \begin{equation}
        \label{errores1}
        \sum \limits_{j=1}^{n} y_{j} = (n+k) \cdot \widehat{\beta}_{1}(k, h) + \widehat{\beta}_{2}(k, h) \cdot \sum \limits_{j=1}^{p} X_{2j} + \cdots  + \widehat{\beta}_{p}(k, h) \cdot \sum \limits_{j=1}^{p} X_{pj} - k \cdot h \cdot \overline{\mathbf{y}}.
    \end{equation}

    On the other hand, based on $\mathbf{e} (k, h) = \mathbf{y} - \mathbf{X} \cdot \widehat{\boldsymbol{\beta}}(k, h)$, it is obtained that:
    $$\mathbf{e}_{j}(k, h) = y_{j} - \left( \widehat{\beta}_{1}(k, h) + \widehat{\beta}_{1}(k, h) X_{2j} + \cdots + \widehat{\beta}_{p}(k, h) X_{pj} \right).$$
   In this case:
    \begin{equation}
        \label{errores2}
        \sum \limits_{j=1}^{n} \mathbf{e}_{j}(k, h) = \sum \limits_{j=1}^{n} y_{j} - \left( n \cdot \widehat{\beta}_{1}(k, h) + \widehat{\beta}_{1}(k, h) \sum \limits_{j=1}^{n} X_{2j} + \cdot + \widehat{\beta}_{p}(k, h) \sum \limits_{j=1}^{n} X_{pj} \right).
    \end{equation}

     Substituting the expression (\ref{errores1}) in (\ref{errores2}), $\sum \limits_{j=1}^{n} \mathbf{e}_{j}(k, h) = k \cdot \widehat{\beta}_{1}(k, h) - k \cdot h \cdot \overline{\mathbf{y}}$,  which in principle differs from zero if $k \not= 0$.

\subsection{Monotony of goodness-of-fit}
    \label{ba.monotonia}

    From expression (\ref{Z.apen}) and taking into account that $\boldsymbol{\Psi} = \mathbf{X}^{t} \mathbf{y}$, it is obtained:
    $$\widehat{\boldsymbol{\beta}}(k, h) = \boldsymbol{\Gamma} \mathbf{D}_{\frac{1}{\lambda_{i}+k}} \boldsymbol{\Gamma}^{t} \left( \boldsymbol{\Psi} + k \cdot h \cdot \boldsymbol{\alpha} \right).$$

    And, then:
    \begin{eqnarray*}
        \widehat{\boldsymbol{\beta}}(k, h)^{t} \left( \mathbf{X}^{t} \mathbf{X} + 2 \cdot k \cdot \mathbf{I} \right) \widehat{\boldsymbol{\beta}}(k, h) &=&  \left( \boldsymbol{\Psi} + k \cdot h \cdot \boldsymbol{\alpha} \right)^{t} \boldsymbol{\Gamma} \mathbf{D}_{\frac{1}{\lambda_{i}+k}} \boldsymbol{\Gamma}^{t} \cdot \left( \boldsymbol{\Gamma} \mathbf{D}_{\frac{1}{\lambda_{i}}} \boldsymbol{\Gamma}^{t} + 2 \cdot k \cdot \boldsymbol{\Gamma} \boldsymbol{\Gamma}^{t} \right)  \\
         & & \boldsymbol{\Gamma} \mathbf{D}_{\frac{1}{\lambda_{i}+k}} \boldsymbol{\Gamma}^{t} \left( \boldsymbol{\Psi} + k \cdot h \cdot \boldsymbol{\alpha} \right) \\
        &=& \mathbf{a}^{t} \mathbf{D}_{\frac{\lambda_{i} + 2 k}{(\lambda_{i}+k)^{2}}} \mathbf{a} + 2 \cdot h \cdot k \cdot \mathbf{a}^{t} \mathbf{D}_{\frac{\lambda_{i} + 2 k}{(\lambda_{i}+k)^{2}}} \mathbf{b} + k^{2} \cdot h^{2} \cdot \mathbf{b}^{t} \mathbf{D}_{\frac{\lambda_{i} + 2 k}{(\lambda_{i}+k)^{2}}} \mathbf{b}. \\
        \widehat{\boldsymbol{\beta}}(k, h)^{t} \boldsymbol{\alpha} &=&  \left( \boldsymbol{\Psi} + k \cdot h \cdot \boldsymbol{\alpha} \right)^{t} \boldsymbol{\Gamma} \mathbf{D}_{\frac{1}{\lambda_{i}+k}} \boldsymbol{\Gamma}^{t} \cdot \boldsymbol{\alpha} \\
        &=& \mathbf{a}^{t} \mathbf{D}_{\frac{1}{\lambda_{i}+k}} \mathbf{b} + k \cdot h \cdot \mathbf{b}^{t} \mathbf{D}_{\frac{1}{\lambda_{i}+k}} \mathbf{b},
    \end{eqnarray*}
    where $\mathbf{a} = \boldsymbol{\Gamma}^{t} \boldsymbol{\Psi}$ and $\mathbf{b} = \boldsymbol{\Gamma}^{t} \boldsymbol{\alpha}$.

    Thus,
    \begin{eqnarray*}
        \widehat{\boldsymbol{\beta}}(k, h)^{t} \left( \mathbf{X}^{t} \mathbf{X} + 2 \cdot k \cdot \mathbf{I} \right) \widehat{\boldsymbol{\beta}}(k, h) - 2 \cdot k \cdot h \cdot \widehat{\boldsymbol{\beta}}(k, h)^{t} \boldsymbol{\alpha} &=& \mathbf{a}^{t} \mathbf{D}_{\frac{\lambda_{i} + 2 k}{(\lambda_{i}+k)^{2}}} \mathbf{a} + 2 \cdot h \cdot \mathbf{a}^{t} \mathbf{D}_{\frac{k^{2}}{(\lambda_{i}+k)^{2}}} \mathbf{b} \\
        & & + h^{2} \cdot \mathbf{b}^{t} \mathbf{D}_{-\frac{\lambda_{i} \cdot k^{2}}{(\lambda_{i}+k)^{2}}} \mathbf{b}.
    \end{eqnarray*}

    Then, expression (\ref{ba.ajuste.alterna}) can be rewritten as:
    \begin{eqnarray}
        GoF(k,h) &=& \frac{\mathbf{a}^{t} \mathbf{D}_{\frac{\lambda_{i} + 2 k}{(\lambda_{i}+k)^{2}}} \mathbf{a} + 2 \cdot h \cdot \mathbf{a}^{t} \mathbf{D}_{\frac{k^{2}}{(\lambda_{i}+k)^{2}}} \mathbf{b} + h^{2} \cdot \mathbf{b}^{t} \mathbf{D}_{-\frac{\lambda_{i} \cdot k^{2}}{(\lambda_{i}+k)^{2}}} \mathbf{b}}{\mathbf{y}^{t} \mathbf{y}} \nonumber \\
        &=& \frac{1}{\mathbf{y}^{t} \mathbf{y}} \cdot \sum \limits_{i=1}^{p} \frac{(\lambda_{i} + 2 k ) a_{i}^{2}}{(\lambda_{i}+k)^{2}} + 2 \cdot h \cdot \frac{1}{\mathbf{y}^{t} \mathbf{y}} \cdot \sum \limits_{i=1}^{p} \frac{k^{2} a_{i} b_{i}}{(\lambda_{i}+k)^{2}} + h^2 \cdot \frac{1}{\mathbf{y}^{t} \mathbf{y}} \cdot \sum \limits_{i=1}^{p} \frac{- \lambda_{i} k^{2} b_{i}^{2}}{(\lambda_{i}+k)^{2}}. \label{ba.sumandos}
    \end{eqnarray}

    Deriving with respect to $k$:
    \begin{eqnarray}
        \frac{\partial}{\partial k} \left( \sum \limits_{i=1}^{p} \frac{(\lambda_{i} + 2 k ) a_{i}^{2}}{(\lambda_{i}+k)^{2}} \right) &=& - \sum \limits_{i=1}^{p} \frac{2 a_{i}^{2} k}{(\lambda_{i}+k)^{3}}, \label{derivada1} \\
        \frac{\partial}{\partial k} \left( \sum \limits_{i=1}^{p} \frac{k^{2} a_{i} b_{i}}{(\lambda_{i}+k)^{2}} \right) &=& \sum \limits_{i=1}^{p} \frac{2 a_{i} b_{i} k \lambda_{i}}{(\lambda_{i}+k)^{3}}, \label{derivada2} \\
        \frac{\partial}{\partial k} \left( \sum \limits_{i=1}^{p} \frac{- \lambda_{i} k^{2} b_{i}^{2}}{(\lambda_{i}+k)^{2}} \right) &=& - \sum \limits_{i=1}^{p}  \frac{2 \lambda_{i}^{2} k b_{i}^{2}}{(\lambda_{i}+k)^{3}}. \label{derivada3}
    \end{eqnarray}

    Due to $k>0$ and $\lambda_{i} > 0$ para $i=1,\dots,p$, it is evident that the derivatives of the expressions (\ref{derivada1}) and (\ref{derivada3}) are negatives, lwhich implies that the first and third summands of expression  (\ref{ba.sumandos}) are decreasing in $k$.
    Note that this implies that for $h=0$ (ridge regression), the goodness-of-fit is decreasing in $k$.

    On the contrary, the sign of the derivative of expression  (\ref{derivada2}) depends on the sign of $a_{i} b_{i}$, for $i=1,\dots,p$, and therefore no statement can be made about the monotonicity of this second summand. . However, it is verified that this term has a horizontal asymptote:
    $$\sum \limits_{i=1}^{p} \frac{k^{2} a_{i} b_{i}}{(\lambda_{i}+k)^{2}} \longrightarrow \sum \limits_{i=1}^{p} a_{i} b_{i} = \mathbf{a}^{t} \mathbf{b} = \boldsymbol{\Psi}^{t} \boldsymbol{\alpha} \mbox{ when } k \rightarrow +\infty,$$
    i.e., there must be a value of $k$ from which the value of this term stabilizes.

    Taking into account the above analysis, it is to be expected that $GoF(k,h)$ decreasing when $k$ increases, aalthough it is not assured. In such case, it would be verified that $2 \cdot h \cdot \boldsymbol{\Psi}^{t} \boldsymbol{\alpha} < GoF(k,h) \leq R^{2}$, where $R^{2}$ is the coefficient of determination of the model (\ref{model0}).

\subsection{Monotony of $\mathbf{S}(k, h)$}
    \label{ba.esgoECM}

    Once again, from the expression (\ref{Z.apen}):
    \begin{eqnarray*}
        \mathbf{Z}(k) \mathbf{X}^{t} \mathbf{X} - \mathbf{I} &=& \boldsymbol{\Gamma} \mathbf{D}_{-\frac{k}{\lambda_{i}+k}} \boldsymbol{\Gamma}^{t}, \\
        \left( \mathbf{Z}(k) \mathbf{X}^{t} \mathbf{X} - \mathbf{I} \right) \mathbf{Z}(k) &=& \boldsymbol{\Gamma} \mathbf{D}_{-\frac{k}{(\lambda_{i}+k)^{2}}} \boldsymbol{\Gamma}^{t}, \\
        \boldsymbol{\beta}^{t} \left( \mathbf{Z}(k) \mathbf{X}^{t} \mathbf{X} - \mathbf{I} \right) \mathbf{Z}(k) \boldsymbol{\alpha} &=& \boldsymbol{\beta}^{t} \boldsymbol{\Gamma} \mathbf{D}_{-\frac{k}{(\lambda_{i}+k)^{2}}} \boldsymbol{\Gamma}^{t} \boldsymbol{\alpha} = \mathbf{e}^{t} \mathbf{D}_{-\frac{k}{(\lambda_{i}+k)^{2}}} \mathbf{b},
    \end{eqnarray*}
    where $\mathbf{e} = \boldsymbol{\Gamma}^{t} \boldsymbol{\beta}$, it can be obtained:
    \begin{equation}
        \label{Shk1}
        2 \cdot k \cdot h \cdot \boldsymbol{\beta}^{t} \left( \mathbf{Z}(k) \mathbf{X}^{t} \mathbf{X} - \mathbf{I} \right) \mathbf{Z}(k) \boldsymbol{\alpha} = - 2 \cdot h \cdot \sum \limits_{i=1}^{p} \frac{e_{i} b_{i} k^{2}}{(\lambda_{i}+k)^{2}}.
    \end{equation}

    On the other hand:
    $$\boldsymbol{\alpha}^{t} \mathbf{Z}(k)^{t} \mathbf{Z}(k) \boldsymbol{\alpha} = \boldsymbol{\alpha}^{t} \boldsymbol{\Gamma} \mathbf{D}_{\frac{1}{(\lambda_{i}+k)^{2}}} \boldsymbol{\Gamma}^{t} \boldsymbol{\alpha} = \mathbf{b}^{t} \mathbf{D}_{\frac{1}{(\lambda_{i}+k)^{2}}}  \mathbf{b},$$
    and, then:
    \begin{equation}
        \label{Shk2}
        k^{2} \cdot h^{2} \cdot \boldsymbol{\alpha}^{t} \mathbf{Z}(k)^{t} \mathbf{Z}(k) \boldsymbol{\alpha} = h^{2} \cdot \sum \limits_{i=1}^{p} \frac{b_{i}^{2} k^{2}}{(\lambda_{i}+k)^{2}}.
    \end{equation}

    From expressions (\ref{Shk1}) and (\ref{Shk2}) it is obtained that:
    $$\mathbf{S}(k, h) = - 2 \cdot h \cdot \sum \limits_{i=1}^{p} \frac{e_{i} b_{i} k^{2}}{(\lambda_{i}+k)^{2}} + h^{2} \cdot \sum \limits_{i=1}^{p} \frac{b_{i}^{2} k^{2}}{(\lambda_{i}+k)^{2}}.$$

    Deriving with respect to $k$:
    \begin{eqnarray}
        \frac{\partial}{\partial k} \left( \sum \limits_{i=1}^{p} \frac{e_{i} b_{i} k^{2}}{(\lambda_{i}+k)^{2}} \right) &=& \sum \limits_{i=1}^{p} \frac{2 k e_{i} b_{i} \lambda{i}}{(\lambda_{i}+k)^{3}}, \label{derivada4} \\
        \frac{\partial}{\partial k} \left( \sum \limits_{i=1}^{p} \frac{b_{i}^{2} k^{2}}{(\lambda_{i}+k)^{2}} \right) &=& \sum \limits_{i=1}^{p} \frac{2 b_{i}^{2} k \lambda_{i}}{(\lambda_{i}+k)^{3}}. \label{derivada5}
    \end{eqnarray}

    Due to $k>0$ and $\lambda_{i} > 0$ for $i=1,\dots,p$, it is clear that the derivative of the expression (\ref{derivada4}) has a sign which depends on $e_{i} b_{i}$, for $i=1,\dots,p$, and the derivative of the expression (\ref{derivada5}) has a positive sign, then the second term of $\mathbf{S}(k, h)$ is increasing in $k$.

    Finally, considering that $k \rightarrow +\infty$:
    $$\mathbf{S}(k, h) \longrightarrow - 2 \cdot h \cdot \sum \limits_{i=1}^{p} e_{i} b_{i} + h^{2} \cdot \sum \limits_{i=1}^{p} b_{i}^{2} = - 2 \cdot h \cdot \mathbf{e}^{t} \mathbf{b} + h^{2} \cdot \mathbf{b}^{t} \mathbf{b} = - 2 \cdot h \cdot \boldsymbol{\beta}^{t} \boldsymbol{\alpha} + h^{2} \cdot \boldsymbol{\alpha}^{t} \boldsymbol{\alpha}.$$
   This is to say, $\mathbf{S}(k, h)$ has a horizontal asymptote, so there must be a value of $k$ at which it stabilizes.

\subsection{Monotony of $|| \alpha - \widehat{\boldsymbol{\beta}} (k,h) ||$}
    \label{rate}

    Since the main objective of the proposed penalty is to ensure that the estimates of $\boldsymbol{\beta}$ are close to $\boldsymbol{\alpha}$, it is interesting to analyze the behavior of $|| \alpha - \widehat{\boldsymbol{\beta}} (k,h) ||$.

    Taking into account thata $\widehat{\boldsymbol{\beta}}(k,h) = \widehat{\boldsymbol{\beta}}(k) + k \cdot h \cdot \mathbf{Z}(k) \boldsymbol{\alpha}$ con $\widehat{\boldsymbol{\beta}}(k) = \mathbf{Z}(k) \mathbf{X}^{t}\mathbf{X} \widehat{\boldsymbol{\beta}}$ (see expression (\ref{penalizado.cresta})), it is obtained that $\boldsymbol{\alpha} - \widehat{\boldsymbol{\beta}}(k,h) = \left( \mathbf{I} - k \cdot h \cdot \mathbf{Z}(k) \right)\boldsymbol{\alpha} - \widehat{\boldsymbol{\beta}}(k)$.
    In this case:
    \begin{eqnarray*}
        || \alpha - \widehat{\boldsymbol{\beta}} (k,h) || &=& \left( \left( \mathbf{I} - k \cdot h \cdot \mathbf{Z}(k) \right)\boldsymbol{\alpha} - \widehat{\boldsymbol{\beta}}(k) \right)^{t}\left( \left( \mathbf{I} - k \cdot h \cdot \mathbf{Z}(k) \right)\boldsymbol{\alpha} - \widehat{\boldsymbol{\beta}}(k) \right) \\
        &=& \boldsymbol{\alpha}^{t} \left( \mathbf{I} - k \cdot h \cdot \mathbf{Z}(k) \right)^{t} \left( \mathbf{I} - k \cdot h \cdot \mathbf{Z}(k) \right)\boldsymbol{\alpha} - 2 \boldsymbol{\alpha}^{t} \left( \mathbf{I} - k \cdot h \cdot \mathbf{Z}(k) \right)^{t} \widehat{\boldsymbol{\beta}}(k) + \widehat{\boldsymbol{\beta}}(k)^{t}\widehat{\boldsymbol{\beta}}(k).
    \end{eqnarray*}

    From $\mathbf{Z}(k) = \boldsymbol{\Gamma} \mathbf{D}_{\frac{1}{\lambda_{i}+k}} \boldsymbol{\Gamma}^{t}$, $\mathbf{X}^{t}\mathbf{X} = \boldsymbol{\Gamma} \mathbf{D}_{\lambda_{i}} \boldsymbol{\Gamma}^{t}$ and remembering that $\mathbf{b} = \boldsymbol{\Gamma}^{t} \boldsymbol{\alpha}$ and $\mathbf{e} = \boldsymbol{\Gamma}^{t} \boldsymbol{\beta}$, it is obtained that:
    \begin{eqnarray*}
        \boldsymbol{\alpha}^{t} \left( \mathbf{I} - k \cdot h \cdot \mathbf{Z}(k) \right)^{t} \left( \mathbf{I} - k \cdot h \cdot \mathbf{Z}(k) \right)\boldsymbol{\alpha} &=& \widehat{\boldsymbol{\alpha}}^{t} \boldsymbol{\Gamma} \mathbf{D}_{\frac{(\lambda_{i} + k - k \cdot h)^{2}}{(\lambda_{i}+k)^{2}}} \boldsymbol{\Gamma}^{t} \widehat{\boldsymbol{\alpha}} = \sum \limits_{i=1}^{p} \frac{(\lambda_{i} + k - k \cdot h)^{2}b_{i}^{2}}{(\lambda_{i}+k)^{2}}, \\
        \boldsymbol{\alpha}^{t} \left( \mathbf{I} - k \cdot h \cdot \mathbf{Z}(k) \right)^{t} \widehat{\boldsymbol{\beta}}(k) &=& \widehat{\boldsymbol{\alpha}}^{t} \boldsymbol{\Gamma} \mathbf{D}_{\frac{\lambda_{i} (\lambda_{i} + k - k \cdot h)}{(\lambda_{i}+k)^{2}}} \boldsymbol{\Gamma}^{t} \widehat{\boldsymbol{\beta}} = \sum \limits_{i=1}^{p} \frac{\lambda_{i} (\lambda_{i} + k - k \cdot h)b_{i}e_{i}}{(\lambda_{i}+k)^{2}}, \\
        \widehat{\boldsymbol{\beta}}(k)^{t}\widehat{\boldsymbol{\beta}}(k) &=& \widehat{\boldsymbol{\beta}}^{t} \boldsymbol{\Gamma} \mathbf{D}_{\frac{\lambda_{i}^{2}}{(\lambda_{i}+k)^{2}}} \boldsymbol{\Gamma}^{t} \widehat{\boldsymbol{\beta}} = \sum \limits_{i=1}^{p} \frac{\lambda_{i}^{2} e_{i}^{2}}{(\lambda_{i}+k)^{2}}.
    \end{eqnarray*}

    Deriving with respect to $k$:
    \begin{eqnarray}
        \frac{\partial}{\partial k} \left( \sum \limits_{i=1}^{p} \frac{(\lambda_{i} + k - k \cdot h)^{2}b_{i}^{2}}{(\lambda_{i}+k)^{2}} \right) &=& - \sum \limits_{i=1}^{p} \frac{2 (\lambda_{i} + k - k \cdot h) b_{i}^{2} h \lambda_{i}}{(\lambda_{i}+k)^{3}}, \label{derivate1} \\
        \frac{\partial}{\partial k} \left( \sum \limits_{i=1}^{p} \frac{\lambda_{i} (\lambda_{i} + k - k \cdot h)b_{i}e_{i}}{(\lambda_{i}+k)^{2}} \right) &=& \sum \limits_{i=1}^{p} \frac{\lambda_{i} (-\lambda_{i} - k - h \lambda_{i} + k \cdot h) b_{i} e_{i}}{(\lambda_{i}+k)^{3}}, \label{derivate2} \\
        \frac{\partial}{\partial k} \left( \sum \limits_{i=1}^{p} \frac{\lambda_{i}^{2} e_{i}^{2}}{(\lambda_{i}+k)^{2}} \right) &=& - \sum \limits_{i=1}^{p} \frac{2 \lambda_{i}^{2} e_{i}^{2}}{(\lambda_{i}+k)^{3}}. \label{derivate3}
    \end{eqnarray}
    From the above expressions it can only be stated that the third term is decreasing in  $k$ (since its derivative, expression (\ref{derivate3}), has a negative sign), so it is difficult to analyze the monotony of the $|| \alpha - \widehat{\boldsymbol{\beta}} (k,h) ||$.
     The one thing that is clear is that when $k \rightarrow +\infty$:
    $$|| \alpha - \widehat{\boldsymbol{\beta}} (k,h) || \rightarrow  (1-h) \cdot \sum \limits_{i=1}^{p} b_{i}^{2} (1-h) \cdot \mathbf{b}^{t} \mathbf{b} = (1-h) \cdot \boldsymbol{\alpha}^{t}\boldsymbol{\alpha} = (1-h) \cdot || \boldsymbol{\alpha} ||.$$

    However, for $h=1$ it is verified that:
    $$|| \alpha - \widehat{\boldsymbol{\beta}} (k,h) || = \sum \limits_{i=1}^{p} \frac{\lambda_{i}^{2}b_{i}^{2}}{(\lambda_{i}+k)^{2}} - 2 \sum \limits_{i=1}^{p} \frac{\lambda_{i}^{2} b_{i}e_{i}}{(\lambda_{i}+k)^{2}} + \sum \limits_{i=1}^{p} \frac{\lambda_{i}^{2} e_{i}^{2}}{(\lambda_{i}+k)^{2}},$$
    which clearly decreases to zero as $k$ increases.

\end{document}